\documentclass[10pt]{article}

\usepackage{latexsym,amssymb,amsmath,theorem,tan-gle,diagrams}
\input diagrams.sty

\newlength{\dinwidth}
\newlength{\dinmargin}
\def\endexem{\hfill{$\Box$}\medskip}

\theoremstyle{change}
\theoremheaderfont{\scshape}
\newtheorem{thm}{Theorem}[section]
\newtheorem{prop}[thm]{Proposition}
\newtheorem{lem}[thm]{Lemma}
\newtheorem{cor}[thm]{Corollary}
\newtheorem{defin}[thm]{Definition}
\newtheorem{defprop}[thm]{Definition/Proposition}
\theorembodyfont{\rmfamily}
\newtheorem{example}[thm]{Example}
\newtheorem{rema}[thm]{Remark}

\newcommand{\bea}{\begin{eqnarray}}
\newcommand{\eea}{\end{eqnarray}}
\newcommand{\bean}{\begin{eqnarray*}}
\newcommand{\eean}{\end{eqnarray*}}

\newcommand{\bdefin}{\begin{defin}}
\newcommand{\blemma}{\begin{lem}}
\newcommand{\bprop}{\begin{prop}}
\newcommand{\btheor}{\begin{thm}}
\newcommand{\bcoro}{\begin{cor}}
\newcommand{\bconj}{\begin{conj}}
\newcommand{\bdefprop}{\begin{defprop}}
\newcommand{\bexam}{\begin{example}}
\newcommand{\edefin}{\end{defin}}
\newcommand{\elemma}{\end{lem}}
\newcommand{\eprop}{\end{prop}}
\newcommand{\etheor}{\end{thm}}
\newcommand{\ecoro}{\end{cor}}
\newcommand{\econj}{\end{conj}}
\newcommand{\brem}{\begin{rema}}
\newcommand{\erem}{\endexem\end{rema}}
\newcommand{\edefprop}{\end{defprop}}
\newcommand{\eexam}{\endexem\end{example}}

\def\1#1{{\bf #1}}
\def\2#1{{\cal #1}}
\def\3#1{{\sl #1}}
\def\4#1{{\tt #1}}
\def\5#1{{\sf #1}}
\def\6#1{{\mathfrak #1}}
\def\7#1{{\mathbb #1}}

\def\qed{\hfill{$\blacksquare$}\medskip}

   \def\d{\delta}

  \def\dd{\partial} \def\D{\Delta}

\def\id{\mathrm{id}}

\newcommand{\End}{\mathrm{End}}
\newcommand{\Hom}{\mathrm{Hom}}

\newcommand{\Rep}{\mathrm{Rep}}
\newcommand{\Mod}{\mathrm{mod}}
\newcommand{\Nat}{\mathrm{Nat}}
\newcommand{\obj}{\mathrm{Obj}}

\newcommand{\Aut}{\mathrm{Aut}}
\newcommand{\mcirc}{\,\circ\,}

\newcommand{\Vect}{\mathrm{Vect}}

\def\prf{\noindent \emph{Proof.\ }}

\newcommand{\restr}{\upharpoonright}
\newcommand{\rarr}{\rightarrow}

\newcommand{\impl}{\Rightarrow}
\newcommand{\ol}{\overline}

\newcommand{\ve}{\varepsilon}
\newcommand{\op}{{\mbox{\scriptsize op}}}

\newcommand{\DS}{\displaystyle}

\newarrow{Congruent} 33333

\numberwithin{equation}{section}

\begin{document}

\title{Monoids, Embedding Functors and Quantum Groups}

\author{Michael M\"uger\footnote{Department of Mathematics, Radboud University, Nijmegen, 
The Netherlands. email: mueger@math.ru.nl}\ \  and Lars Tuset\footnote{Faculty of Engineering, Oslo University College, Oslo, Norway. email: Lars.Tuset@iu.hio.no}}

\date{\today}
\maketitle

\begin{abstract}
We show that the left regular representation $\pi_l$ of a discrete quantum group 
$(A,\D)$ has the absorbing property and forms a monoid $(\pi_l,\tilde{m},\tilde{\eta})$ 
in the representation category $\Rep(A,\D)$.

Next we show that an absorbing monoid in an abstract tensor $*$-category $\2C$ gives rise to an 
embedding  functor (or fiber functor) $E:\2C\rarr\mathrm{Vect}_\7C$, and we identify conditions 
on the monoid, satisfied by $(\pi_l,\tilde{m},\tilde{\eta})$, implying that $E$ is $*$-preserving.  

As is well-known, from an embedding functor $E: \2C\rarr\mathrm{Hilb}$ the 
generalized Tannaka theorem produces a discrete quantum group $(A,\D)$ 
such that $\2C\simeq\Rep_f(A,\D)$.
Thus, for a $C^*$-tensor category $\2C$ with conjugates and irreducible unit the following are equivalent:
(1) $\2C$ is equivalent to the representation category of a discrete quantum 
group $(A,\D)$, 
(2) $\2C$ admits an absorbing monoid,
(3) there exists a $*$-preserving embedding functor $E: \2C\rarr\mathrm{Hilb}$.
\end{abstract}


\section{Introduction and related work}
\subsection{Our approach}\label{ss-int}
As is well-known, see for example \cite[Sections 2-3]{MRT},
the finite dimensional representations of a discrete quantum group 
form a $C^*$-tensor category with conjugates and irreducible unit.  
It is therefore natural to ask for a characterization of representation categories of 
discrete quantum groups among the $C^*$-tensor categories.
A partial solution is provided by the generalized Tannaka theorem, cf.\ \cite{woro, JS}, according
to which a $C^*$-tensor category is such a representation category whenever it comes
equipped with an embedding functor, i.e.\ a faithful $*$-preserving tensor functor into the category
$\2H$ of finite dimensional Hilbert spaces.
In this case the category is called concrete as opposed to abstract. 
The most transparent approach to the Tannaka theorem defines the 
quantum group as the algebra of natural transformations of the embedding
functor to itself. The monoidal structures of the category and of the embedding functor then give
rise to the coproduct of the  quantum group. For this approach and further references cf.\  \cite{MRT}.  

The generalized Tannaka theorem reduces the characterization problem to that of producing 
an embedding functor. Since the representation category of a quantum group
comes with an obvious embedding functor, the existence of such a functor
clearly is a necessary condition. However, there exist $C^*$-tensor 
categories with conjugates and irreducible unit that do not admit an 
embedding functor: Infinitely many examples (which are even braided) are 
provided by the categories associated with quantum groups at roots of unity, 
cf.\ \cite{Wenzl}. This shows that additional assumptions on an abstract $C^*$-tensor category are
needed in order to identify it as 
the representation category of a quantum group. For example, in \cite{Wenzl2} it is proven that
any $C^*$-tensor category with conjugates, irreducible unit and with fusion ring
isomorphic to that of $SU(N)$ is equivalent to the representation category
of the discrete quantum group dual to $SU_q(N)$ for some $q\in\7R$. Analogous results
have been proven for the other classical groups, assuming in addition
that the category is braided. 

The case of abstract symmetric tensor categories was settled already in the late 80's.
By a remarkable result of Doplicher and Roberts \cite{DR}, any symmetric
$C^*$-tensor category with conjugates and irreducible unit is equivalent
as a $C^*$-tensor category to the representation category of a unique compact group.
If one wishes an equivalence of symmetric categories, one must also
allow super groups. This result has applications \cite{DR2} to algebraic quantum field theory,
where symmetric $C^*$-tensor categories arise without an a priori given embedding functor. 
The proof in \cite{DR}, however, does not follow the strategy outlined above of constructing an
embedding functor and then applying the Tannaka theorem. 

Independently and at about the same time, motivated by applications to algebraic geometry,
Deligne proved \cite{del} that a rigid abelian symmetric tensor category with irreducible unit 
is equivalent to the representation category of a proalgebraic group, provided that the intrinsic
dimension of every object is a positive integer. His proof consists of constructing an embedding
functor and applying the algebraic Tannaka  theorem of N. Saavedra Rivano. 

The crucial notion in Deligne's construction of the embedding functor is that of an absorbing
commutative monoid. Recall that a monoid in a tensor category is a triple $(Q,m,\eta)$, where 
$m:Q\otimes Q\to Q$ and $\eta :\11\to Q$ are morphisms such that 
$(m\otimes\id_Q )\circ m=(\id_Q\otimes m)\circ m$ and 
$m\circ (\eta\otimes\id_Q)=\id_Q =m\circ (\id_Q\otimes\eta)$. An object $Q$ 
is called absorbing if the $Q$-module $Q\otimes X$ is isomorphic to some multiple of $Q$ for 
any object $X$. Deligne obtained the absorbing commutative monoid using categorical 
generalizations of results from commutative algebra -- it is here that the symmetry 
plays a central role. His proof was simplified considerably in \cite{bichon}. Note, however, that
the monoid of \cite{del,bichon} fails to satisfy hypothesis 1 of Proposition \ref{p-embed} below, which
complicates the construction of an embedding functor. For a construction of a monoid satisfying all
assumptions of Proposition \ref{p-embed} cf.\ \cite{MM}.  
   
\medskip

The aim of this paper is to demonstrate the usefulness of the monoid approach
in the general non-symmetric case. This is done in two steps. On the one hand
we prove that the passage from an absorbing monoid to an embedding functor works in
the general case. We also identify conditions on the monoid guaranteeing that the functor
is $*$-preserving. 
Whereas the existence of an embedding functor refers to $\2H$ and thus is an external 
condition on the category, the existence of an absorbing monoid is an internal property.
As such it is more amenable to proof, as Deligne's result in the symmetric case illustrates.
A technical aspect should be pointed out though: A category $\2C$ with conjugates can contain
an absorbing object only if it has finitely many equivalence classes of objects. 
Otherwise it needs to be suitably enlarged, which is done using the category $\hat{\2C}$ 
of inductive limits. We say that $\2C$ {\it admits} an absorbing object if there exists a monoid
$(Q,m,\eta)$ in $\hat{\2C}$ such that the $Q$-module $(Q\otimes X,m\otimes\id_X)$ is isomorphic to a
multiple of the $Q$-module $(Q,m)$, for every $X\in\2C$.

On the other hand, starting with a discrete quantum group $(A,\D)$, we explicitly construct an
absorbing monoid $(\pi_l,\tilde{m},\tilde{\eta})$ in the representation category. Here
$\pi_l$ is the regular representation of the algebra $A$ on the vector space $A$ given by
multiplication from the left. In order to define the morphisms $\tilde{m},\tilde{\eta}$, 
let $(\hat{A}, \hat{\D})$ denote the dual compact quantum group with multiplication $\hat{m}$ and
unit $1_{\hat{A}}$, and let $\2F :A\to\hat{A},\ a\mapsto\varphi(\cdot\,a)$ denote the Fourier
transform, where $\varphi :A\to\7C$ is the left invariant positive functional of $(A,\D)$. 
The linear maps $\tilde{m}:A\otimes A\to A$ and $\tilde{\eta}:\7C\to A$ are then given by
$\tilde{m}=\2F^{-1}\hat{m}(\2F\otimes\2F)$ and $\tilde{\eta}(1)=\2F^{-1}(1_{\hat{A}})$.
We call this absorbing monoid the regular monoid of $(A,\D)$.

Our main result then is that, for a $C^*$-tensor category $\2C$ with conjugates and 
irreducible unit, we have three equivalent statements illustrated by the 
following diagram:

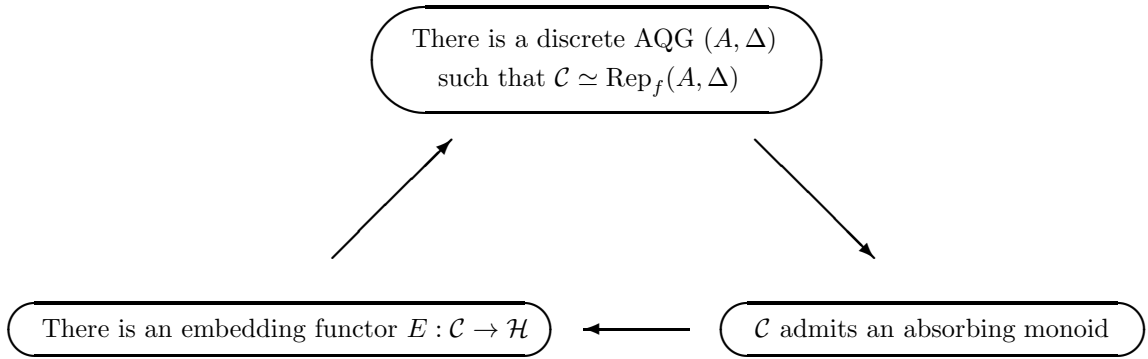
\[\begin{array}{c}\begin{picture}(300,150)(80,0)\thicklines

\put(160,130){There is a discrete AQG $(A,\D)$}

\put(170,115){such that $\2C\simeq\Rep_f(A,\D)$}

\put(230,125){\oval(170,40)}

\put(290,20){$\2C$ admits an absorbing monoid}

\put(357,23){\oval(160,20)}

\put(20,20){There is an embedding functor $E:\2C\rarr\2H$}

\put(110,23){\oval(205,20)}

\put(130,50){\vector(1,1){45}}

\put(290,95){\vector(1,-1){45}}

\put(265,23){\vector(-1,0){40}}

\end{picture} \end{array}\]

We summarize some further results. Our construction actually provides an absorbing semigroup 
$(\pi_l , \tilde{m})$ for any algebraic quantum group, and we show that this semigroup has a unit 
$\tilde{\eta}$ if and only if the quantum group is discrete. Dually, there exists a regular 
comonoid if and only if the quantum group is compact. 
In the finite dimensional case the regular monoid and comonoid combine to a Frobenius algebra. 
We identify the intrinsic group of a discrete quantum group with the intrinsic group of its
regular monoid. 
     
We also show that an abstract $C^*$-tensor category $\2C$ with conjugates and irreducible unit
admits an absorbing object $Q$ in $\hat{\2C}$ if and only if $\2C$ admits an integer valued 
dimension function, i.e.\ a map $\obj\,\2C\to\7N$ that is additive and multiplicative. While this
clearly is a necessary condition for $\2C$ to admit an absorbing monoid,  
to proceed further in the opposite direction one also needs an associative 
morphism $m:Q\otimes Q\to Q$, but the existence of such a morphism remains to be proven. 

\subsection{Related work}
We would like to point out several earlier references that are related to the present work. The fact that 
a $C^*$-tensor category with finitely many simple objects and an absorbing monoid is the representation 
category of a finite dimensional $C^*$-Hopf algebra was obtained in \cite[Theorem 6.7]{LR}. The proof basically 
proceeds by showing that a finite $C^*$-tensor category can be faithfully realized by endomorphisms of a von 
Neumann algebra and then appealing to \cite[Theorem 6.2]{longo}. This approach can in principle (this has never 
been done) be extended to $C^*$-tensor categories with countably many simple objects by combining \cite{yama2}, 
which realizes such categories as categories of bimodules over a von Neumann algebra $N$, and the extension of 
Longo's result \cite{longo} to infinite index subfactors of depth two obtained in \cite{FI}. Concerning this
generalization we observe that putting the above-mentioned results together as indicated would require a
non-trivial amount of work since they use different frameworks (type II vs.\ III algebras, endomorphisms 
vs.\ bimodules), and dropping the countability assumption on the category made in \cite{yama2} seems very 
difficult. Furthermore, the above approaches (in the finite and countable cases) use very heavy operator algebraic
machinery, whereas the approach outlined in Subsection \ref{ss-int} is essentially purely algebraic and quite 
elementary and has the added benefit of working without any assumption on the cardinality of the category.

More recently, the relationship between the absorbing property and embedding functors has been  studied in
\cite[Appendices A-B]{yama1}, though with different emphasis and results. Finally, we'd like to point out the
papers \cite{dpr}, which provide a study, in the context of $C^*$-tensor categories, of multiplicative unitaries,  
which are a convenient tool for the study of the regular representation in the theory of locally 
compact groups and quantum groups.


\section{From Algebraic Quantum Groups to Absorbing Monoids} \label{s-regular}
\subsection{Three Representation Categories of AQG}
For the general theory of algebraic quantum groups (AQG) 
we refer to \cite{VD} and to, e.g., \cite{ddz} for the basics of
representation theory, as well as to the survey \cite{MRT}, where both subjects are covered in
considerable detail. For the standard categorical notions of (braided/symmetric) tensor categories,
natural transformations etc., our standing reference is \cite{cwm}, but most of the relevant notions
can also be found in \cite{MRT}. We will always denote AQG by $(A,\D)$, where $A$ is a
non-degenerate $*$-algebra and $\D: A\rarr M(A\otimes A)$ is the comultiplication. 
As usual, we denote the multiplication, counit and coinverse by 
$m: A\otimes A\rarr A,\ \ve: A\rarr\7C$ and $S: A\rarr A$, respectively.
The left invariant positive functional is denoted by $\varphi$.

\bdefin \label{d-repn}
Let $(A, \D )$ be an AQG. A homomorphism $\pi: A\rarr \End\,K$, where $K$ is a complex vector space,
is called a representation of $A$ on $K$ if $\pi(A)K=K$. A $*$-representation is a representation
$\pi$ on a pre-Hilbert space $K$, that is $*$-preserving in the sense that
$(\pi(a)u,v)=(u,\pi(a^*)v)$ for all $a\in A$ and $u,v\in K$. By $\Rep\, (A, \D )$ we denote the
category whose objects are $*$-representations and whose arrows are the intertwining linear maps,
i.e.\ if $\pi '$ is another $*$-representation of $A$ on $K'$, then
\[ \Hom(\pi ,\pi' )=\{s\in \Hom(K,K') \ | \ s\pi(a)v=\pi'(a)sv\ \forall v\in K, \ a\in A\}. \]
\edefin

\brem 
Recall that a homomorphism $\pi: A\rarr B$ of non-degenerate algebras is called non-degenerate if
$\pi(A)B=B=B\pi(A)$. It would therefore seem natural to define a representation of $A$ on $K$ to be
a homomorphism $\pi: A\rarr\End K$ that satisfies $\pi(A)\End K=\End K=\End K\pi(A)$. However, this
notion is too restrictive since it is never satisfied by the usual left regular representation
$\pi_l$, to be introduced shortly, if $(A,\D)$ is discrete and non-unital. So see this it suffices
to notice that the image of $\pi_l(a)e\in\End A$ is finite dimensional for all $a\in A$ and
$e\in\End A$.
\erem

We define the left regular representation $\pi_l: A\rarr\End\,A$ of an AQG by
$\pi_l(a)(x)=ax$ for $a,x\in A$. This terminology is justified, since the non-degeneracy
condition in Definition \ref{d-repn} holds because $A^2=A$, which again follows from
the existence of local units for $A$. Furthermore, $\pi_l$ is a $*$-representation with respect to
the inner product $(\cdot ,\cdot )$ on $A$ given by $(x,y)=\varphi (y^* x)$.
Thus $\pi_l\in\Rep(A,\D)$. Similarly, one defines the right regular representation
$\pi_r\in\Rep(A_\op,\D)$ by the formula $\pi_r (a)(x) =xa$ for $a,x\in A$. 
It is a $*$-representation with respect to the inner product on the opposite algebra 
$A_\op$ given by $(x,y)=\varphi (xy^*)$.

Recall that the left multiplier algebra $L(A)$ of a non-degenerate algebra $A$ is the vector space 
$L(A)=\{\psi\in\End\,A \ |\ \psi (ab)=\psi (a)b \ \forall a,b\in A\}$ with product
$\psi_1\psi_2=\psi_1\circ \psi_2$, i.e.\ composition of maps. Note that $\pi_l: A\rarr L(A)$ is an
injective algebra homomorphism. Similarly, the right multiplier algebra $R(A)$ of a non-degenerate
algebra $A$ is the vector space 
$R(A)=\{\phi\in\End\,A\ |\ \phi (ab)=a\phi (b) \ \forall a,b\in A\}$ together with
the product given by opposite composition: $\phi_1\phi_2=\phi_2\circ \phi_1$. 
Again $\pi_r: A\rarr R(A)$ is an injective algebra homomorphism. Further, note that by definition  
the identity map is a linear antimultiplicative map from $\Hom(\pi_l,\pi_l)$ to $R(A)$ and a linear
multiplicative map from $\Hom(\pi_r,\pi_r)$ to $L(A)$. Assume $A$ is a non-degenerate $*$-algebra
and let $\psi\in L(A)$. Define $\psi^*\in R(A)$ by $\psi^* (a)=\psi (a^*)^*$ for $a\in A$. The
assignment $\psi\mapsto\psi^*$ is a antilinear and antimultiplicative bijection from $L(A)$ to
$R(A)$. The multiplier algebra $M(A)$ of a non-degenerate algebra $A$ is the vector space
$M(A)=\{(\psi ,\phi )\in L(A)\times R(A)\ |\ \phi (a)b=a\psi (b) \ \forall a,b\in A\}$ with
pointwise multiplication, i.e.\ 
$(\psi_1,\phi_1)(\psi_2,\phi_2)=(\psi_1\psi_2,\phi_1\phi_2)=(\psi_1\circ\psi_2,\phi_2\circ\phi_1)$.
Now the map $\pi_{lr} :a\mapsto (\pi_l (a),\pi_r (a))$ embeds $A$ into $M(A)$ as an algebra.
Whenever $A$ is a $*$-algebra, so is $M(A)$ and the embedding is $*$-preserving. If $A$ is unital 
then we have the algebra isomorphisms $M(A)\cong L(A)\cong R(A)\cong A$.

Any homomorphism $\pi: A\rarr\End K$ of a non-degenerate algebra $A$ such that $\pi (A)K=K$ and 
such that $\pi (A)v =0$ implies $v=0$ has a unique extension to a unital homomorphism 
$\tilde{\pi}:M(A)\rightarrow \End K$ given by the formula
$\tilde{\pi} (x)\pi (a)v =\pi (xa)v$, for $x\in M(A)$, $a\in A$ and $v\in K$. Whenever $A$
has local units, the property $\pi (A)v =0\ \impl\ v=0$ follows immediately from $\pi (A)K=K$,  
see \cite{ddz} for more details. If $\pi ,\pi'\in\Rep(A, \D )$, then clearly  
$\pi\otimes\pi' :A\otimes A\rarr\End K\otimes\End K'\subset\End(K\otimes K')$ determined by
$(\pi\otimes\pi')(a\otimes a' )=\pi (a)\otimes\pi (a' )$ for $a,a'\in A$ satisfies 
$(\pi\otimes\pi')(A\otimes A)(K\otimes K')=K\otimes K'$. It therefore has a unique
extension to a unital $*$-homomorphism from $M(A\otimes A )$ to $\End(K\otimes K')$, which
we again denote by $\pi\otimes\pi'$. It is obvious that $\pi\times\pi'=(\pi\otimes\pi')\circ\D$ is
non-degenerate, and therefore belongs to $\Rep(A,\D)$. Hence $\Rep(A, \D )$ is a tensor
category with irreducible unit $\varepsilon$. 
Suppressing the totally canonical associativity constraint, we treat
the tensor category $\Rep(A,\D)$ as strict.
Note that $(\pi_l\times\pi_l )(a)x =\Delta (a)x$ for $a\in A$ and 
$x\in A\otimes A$. By $\Rep_f (A,\D)$ we mean the full tensor subcategory of $\Rep\, (A,\D)$
consisting of finite dimensional representations, i.e. those $\pi\in\Rep\, (A,\D)$ for which 
$\dim K <\infty$.  

Clearly, $\Rep_f (A,\D)$ is a tensor $*$-category w.r.t.\ the adjoint operation for bounded linear
maps between Hilbert spaces, but we are not aware of a method to turn $\Rep\, (A, \D )$
into a tensor $*$-category which works for any AQG $(A,\D)$. Yet, we have the following.

\bprop \label{p-rep*}
Let $(A,\D )$ be an AQG and define $\Rep_*(A,\D)$ to be the full subcategory of $\Rep(A,\D)$ 
consisting of representations that are direct sums of finite dimensional irreducible
$*$-representations with finite multiplicities. Then there exists a $*$-operation on $\Rep_*(A,\D)$
extending that of $\Rep_f(A,\D)$. This $*$-operation is compatible with the scalar products in the
sense that 
\[ (su,v)_{K'}=(u,s^*v)_K \]
for $u\in K, v\in K'$ and $s\in\Hom(\pi,\pi')$, where $\pi,\,\pi'$ are representations on 
$K, \,K'$ with inner products $(\cdot,\cdot)_K$ and $(\cdot,\cdot)_{K'}$, respectively. For
$\pi\cong\oplus_i n_i \pi_i$ and $\pi'\cong\oplus_i n'_i \pi_i$, where the
representations $\pi_i\in\Rep_f(A,\D)$ are irreducible and pairwise non-isomorphic, we use the isomorphisms 
\[ \Hom(\pi,\pi')\cong\prod_i \Hom(n_i\pi_i,n'_i\pi_i)\cong\prod_i M_{n_i,n'_i}(\7C) \]
to equip the spaces $\Hom(\pi,\pi')$, where $\pi,\pi'\in\Rep_*(A,\D)$, with the product
topology. With respect to these topologies the composition $\circ$ is continuous.
\eprop

\prf Let $I$ be the set of unitary equivalence classes of finite dimensional irreducible
$*$-representations and let $\pi_i$ be a representation in the class $i\in I$ acting on the Hilbert
space $H_i$. Consider two representations 
$\pi\cong\oplus_i \pi_i\otimes I_{K_i}$ and $\pi'\cong\oplus_i \pi_i\otimes I_{K'_i}$
where $K_i,K'_i$ are finite dimensional multiplicity spaces. Here it is understood that the scalar
products on the finite dimensional spaces $H_i\otimes K_i$ are the restrictions of that of $K$ and
similarly for $K'$, etc. Since the representations $\pi_i\otimes I_{K_i}$ and 
$\pi_j\otimes I_{K'_j}$ are disjoint if $i\ne j$, every morphism $s: \pi\rarr\pi'$ is given by a
family $(s_i)$, where $s_i\in\Hom(\pi_i\otimes I_{K_i},\pi_i\otimes I_{K'_i})$. Here $s_i$ is a
morphism in the $*$-category $\Rep_f(A,\D)$ and therefore has an adjoint $s_i^*$ defined by
$(s_iu_i,v_i)_{H_i\otimes K'_i}=(u_i,s_i^*v_i)_{H_i\otimes K_i}$. Conversely, every such family
constitutes a morphism in $\Hom(\pi,\pi')$. Thus we can define an element of $\Hom(\pi',\pi)$ by
$s^*=(s_i^*)$. It is evident that this definition satisfies the
properties of a $*$-operation and extends the $*$-operation of $\Rep_f(A,\D)$. Now
$(su,v)_{K'}=(u,s^*v)_K$ is automatic since $(\cdot,\cdot)_K=\sum_i(\cdot,\cdot)_{H_i\otimes K_i}$,
etc. The continuity of $\circ$ is also obvious.
\qed

\brem 1. Note that $\Rep_*(A,\D)$ is not closed under tensor products, but it is
stable under tensor products with finite dimensional $*$-representations.

2. For a general AQG the category $\Rep_*(A,\D)$ may consist only of copies of $\varepsilon$. 
This does not happen in the discrete case to be discussed below.
\erem

\bprop \label{p-disc}
Let $(A,\D)$ be a discrete AQG, so $A=\oplus_{i\in I} \End\, H_i$ with $H_i$ finite
dimensional Hilbert spaces. Let $I_i$ be the unit of $\End\, H_i$ and let $p_i\in\Rep\, (A,\D)$ denote
the canonical projection from $A$ to $\End\, H_i$. Then: 
\begin{enumerate}
\item For any $\pi\in\Rep\, (A,\D)$ we have $\pi\cong\oplus_i n_i p_i$ with 
$n_i =\dim \pi (I_i )K /\dim H_i$.
\item $\DS \pi_l\cong\oplus_{i\in I} \dim H_i\ p_i$, so $\pi_l\in\Rep_*(A,\D)$. 
\item $\Rep_f (A,\D)$ is equivalent to the tensor category of all finite dimensional representations
of $(A,\D)$ and $\Rep(A,\D)$ is equivalent to the tensor category of all representations. 
\item $R(A)\cong M(A)$ as unital algebras, whereas $\Hom(\pi_l,\pi_l)$ and $M(A)$ are 
anti-isomorphic as unital $*$-algebras.
\end{enumerate}
\eprop

\prf 
1. The subspaces $K_i =\pi (I_i )K$ are clearly linear independent and a short argument using $\pi (A)K=K$
shows that $K\cong\oplus_i K_i$. Define $*$-representations $\pi_i$ of $(A,\D )$ on $K_i$
by $\pi_i (a)=\pi (a)\restr K_i$ for $a\in A$, and note that $\pi\cong\oplus_i\pi_i$ with 
$\pi_i\cong n_i p_i$. 

2. This follows from 1. by noting that $K_i =\End\, H_i$ so $n_i =\dim H_i$.

3. This follows from the facts that the decomposition in 1. holds also for representations which are not
$*$-representations and that the irreducible representations $p_i$ are $*$-representations.

4. By definition $\Hom(\pi_l,\pi_l)$ and $R(A)$ are anti-isomorphic as unital algebras. 
Let $\phi\in R(A)$. In view of the definition of
right multipliers we have $\phi(a)=\phi(I_ia)=I_i\phi(a)\in \End\, H_i$ for every 
$i\in I, \ a\in \End\, H_i$. Thus we obtain restrictions 
$\phi_i=\phi\restr \End\, H_i\in R(\End\, H_i)$ such that $\phi=\oplus_i \phi_i$. Conversely,
the latter formula defines an element of $R(A)$ for every element $(\phi_i, i\in I)$ of 
$\prod_i R(\End\, H_i)$. Since the $\End\, H_i$ are unital we have $R(\End\, H_i)=\End\, H_i$, and
therefore $R(A)=\prod_i \End\, H_i=M(A)$ as unital algebras. It follows that $\Hom(\pi_l,\pi_l)$ and
$M(A)$ are anti-isomorphic as unital $*$-algebras. 
\qed

For a discrete AQG we normalize the left invariant positive functional by requiring
$\varphi(I_0)=1$. 

We aim now at understanding the relation between  $\Rep_f(A,\D)$ and $\Rep(A,\D)$
in more categorical terms, whenever $(A,\D)$ is a discrete quantum group.

In order to make sense of infinite direct sums of objects we need some categorical devices.
Let $\2J$ be a small category, the index category, and let $F: \2J\rarr\2C$ be a functor. We denote
the objects of $\2J$ by $i,j,k$ and write $X_i=F(i)$. A pair $(X,f_i)$, where $X\in\2C$ and the
morphisms $f_i: X_i\rarr X$ for $i\in\2J$ satisfy $f_j\circ F(s)=f_i$ for every $s: i\rarr j$, is
called a cone. We say $F$ has an {\it inductive limit} (or colimit) if there exists a cone
$(X,f_i)$ that is universal, i.e., for any other cone $(Y,g_i)$ there exists a unique
$t: X\rarr Y$ such that $t\circ f_i=g_i$ for all $i\in\2J$. 
The category $\2J$ is {\it filtered} if it satisfies the following conditions: 
\begin{enumerate} 
\item For every $i,j\in\2J$, there exists $k\in\2J$ and morphisms $u: i\rarr k$ and $v: j\rarr k$.
\item For every $i,j\in\2J$ and $u,v: i\rarr j$, there exists $s: j\rarr k$ such that
$s\circ u=s\circ v$. 
\end{enumerate} 
An inductive limit $F: \2J\rarr\2C$ is called {\it filtered} if $\2J$ is a filtered category.
Every directed partially ordered set $J$ gives rise to a filtered category $\2J$, where 
$\obj\,\2J=J$ and $\Hom_\2J(i,j)$ contains one element if $i\le j$ and none otherwise.
Given a set $S$, the power set $2^S$ is a directed partially ordered set.

In our applications $\2C$ has finite direct sums, and we define an infinite direct sum
$\oplus_{j\in S} Y_j$ as a filtered inductive limit over $F: 2^S\rarr\2C$. Here $2^S$ is the
filtered category corresponding to the power set $2^S$, and the functor $F$ is
given by choosing a finite direct sum for every $s\in 2^S$. An example of a category for which all
filtered inductive limits exist is the category $\Rep(A,\D)$, where $(A,\D )$ is a discrete AQG. 

We will now consider a completion $\hat{\2C}$ w.r.t.\ all filtered inductive limits of a given
category $\2C$. Given any category $\2C$ there exists a category $\mathrm{Ind}\,\2C$ of `filtered
inductive limits of objects in $\2C$'. The standard reference is \cite{SGA4}. We collect some of its
properties that we shall need, none of which is new. 

\bprop \label{p-ind}
Let $\2C$ be a category and denote $\hat{\2C}=\mathrm{Ind}\,\2C$. Then
\begin{enumerate}
\item $\hat{\2C}$ contains $\2C$ as a full subcategory.
\item $\hat{\2C}$ is complete w.r.t.\ filtered inductive limits. In particular, there exist infinite
sums $Z\cong\bigoplus_{i\in I} Z_i$, where $Z_i\in\2C$.
\item If $\2C$ is abelian, in particular semisimple, then $\hat{\2C}$ is abelian.
\item If $\2C$ is semisimple then every object of $\hat{\2C}$ is a filtered inductive limit of
objects in $\2C$. In this case, $\hat{\2C}$ is uniquely characterized up to equivalence by this
property and 1-2.
\item If $\2C$ is monoidal then the tensor product extends uniquely to $\hat{\2C}$. Similarly if
$\2C$ is braided or symmetric, then so is $\hat{\2C}$.
\item If $\2C$ has exact tensor product, in particular if $\2C$ has duals, then the tensor product
of $\hat{\2C}$ is also exact.
\end{enumerate}
\eprop

\prf We limit ourselves giving references for the interested reader. 
Statements 1-2 are proven in \cite{SGA4}, whereas 3. follows from \cite{SGA4,gab}. 
Claim 4 is proven in \cite[\S 4]{del2}, and for 5-6. see \cite{DM,del}. 
\qed

\brem Concerning the construction of $\hat{\2C}$ we only note that its objects are pairs $(\2G, F)$, 
where $\2G$ is a small filtered category and $F: \2G\rarr\2C$ is a functor. Denoting
objects of $\hat{\2C}$ by $(X_i)$, where $i\in\obj\,\2G$ and $X_i=F(i)$, the hom-sets are defined by 
\[ \Hom_{\hat{\2C}}((X_i),(Y_j))=\varprojlim_i\varinjlim_j\, \Hom_\2C(X_i,Y_j). \]
\erem

\blemma Let $\2C$ be a semisimple tensor $*$-category, and let $\2C_*$ be the full subcategory of
$\hat{\2C}$ consisting of direct sums of irreducible objects of $\2C$ with finite
multiplicities. Then $\2C_*$ has a $*$-operation extending that of $\2C$.
\elemma

\prf Exactly as for Proposition \ref{p-rep*}.
\qed

The following is an immediate, though very useful generalization of \cite[Example 4.3.2]{del2}.

\bprop 
Let $(A,\D)$ be a discrete AQG and let $\2C=\Rep_f(A,\D)$. Then there is a canonical equivalence
$\hat{F}: \hat{\2C}\rarr\Rep(A,\D)$ of tensor categories which restricts to the identity on the full 
subcategory $\2C$ of $\hat{\2C}$ and restricts to an equivalence of $\2C_*$ and $\Rep_*(A,\D)$.
\eprop

\prf Note that if $\2C$ is a category of vector spaces or of representations, then the filtered
inductive limits above are inductive limits in the ordinary sense.
The category $\Rep_f(A,\D)$ is semisimple and every object of
$\Rep(A,\D)$ is an inductive limit of objects in $\Rep_f(A,\D)$. Since $\Rep(A,\D)$ is closed
w.r.t.\ inductive limits, the equivalence $\hat{\2C}\simeq\Rep(A,\D)$ follows from assertion 4
in Proposition \ref{p-ind}. The last statement is obvious since both $\2C_*$ and $\Rep_*(A,\D)$ are
defined as the respective full subcategories of objects that contain the simple objects with finite
multiplicities.
\qed


\subsection{Construction of the Regular Monoid}
Let $(A,\D )$ be an AQG and $(\hat{A},\hat{\D})$ its Pontryagin dual with the conventions that
$\hat{\Delta}(\omega)(a\otimes b)=\omega(ab)$ and 
$\hat{m}(\omega\otimes\omega')=(\omega\otimes\omega')\Delta$, where $a,b\in A$ and
$\omega,\omega'\in\hat{A}$.
Consider the Fourier transform $\2F :A\rarr\hat{A}$, which is given by 
$\2F (a)=\hat{a}=a\varphi$, for $a\in A$. 
Here and in the sequel $c\varphi$ and $\varphi c$ denote the linear functionals on $A$ given by
$c\varphi =\varphi (\cdot\, c)$ and $\varphi c=\varphi (c\, \cdot)$, for $c\in M(A)$.
It is known that $\2F$ is a
bijective linear map satisfying Plancherel's formula
$\hat{\psi}(\2F(a)^* \2F(b))=\varphi (a^* b)$, for $a,b\in A$.
Here $\hat{\psi}$ is the right invariant functional on $(\hat{A} ,\hat{\D})$ determined by
$\hat{\psi}\2F=\varepsilon$.      
If $(A,\D )$ is discrete, then $(\hat{A},\hat{\D})$ is a Hopf $*$-algebra and $\hat{\psi}$
is a bounded functional on $\hat{A}$ which is both left- and right invariant.

\blemma
\label{l-strinv}
Let $(A,\D)$ be an AQG. Then
\[(\varphi\otimes\varphi )(\D (c)(a\otimes b))
=\varphi (c (\varphi\otimes\iota )[((S^{-1}\otimes\iota )\D (b))(a\otimes 1)])
=\varphi (c(\iota\otimes\varphi )[(1\otimes S^{-1}(b))\D (a)]),\]
for $a,b\in A$ and $c\in M(A)$.
\elemma
\prf
The formula $\varphi ((\omega S\otimes\iota)\D (c)b)=\varphi (c (\omega\otimes\iota )\D (b))$ holds for
any $\omega\in\hat{A}$, $b\in A$ and $c\in M(A)$, and is known as the strong left 
invariance property \cite{KD}. Thus
\begin{eqnarray*} \lefteqn{
(\varphi\otimes\varphi )(\D (c)(a\otimes b))
=\varphi(((\hat{a}S^{-1} S\otimes\iota )\D (c))b) } \\
 && =\varphi (c(\hat{a}S^{-1} \otimes\iota )\D (b))
=\varphi (c (\varphi\otimes\iota )[((S^{-1}\otimes\iota )\D (b))(a\otimes 1)])
\end{eqnarray*}
for $a,b\in A$ and $c\in M(A)$. The computation
\begin{eqnarray*}  \lefteqn{
   \varphi (c (\varphi\otimes\iota )[((S^{-1}\otimes\iota )\D (b))(a\otimes 1)])
  = \varphi (c(\hat{a}S^{-1} \otimes\iota )\D (b)) } \\
 && = \hat{a}S^{-1}((\iota\otimes\varphi c)\D (b)) =
   \varphi (S^{-1}[(\iota\otimes\varphi c)\D (b)]a) \\
 && = \varphi ((\varphi cS\otimes\iota )\D S^{-1}(b) a) =
   \varphi (S^{-1}(b)(\varphi c\otimes\iota)\D (a)) \\
 && =\varphi (c(\iota\otimes\varphi )[(1\otimes S^{-1}(b))\D (a)])
\end{eqnarray*}
proves the second identity.
\qed

Let $\hat{m} :\hat{A}\otimes\hat{A}\rarr\hat{A}$ be the linearized multiplication on $\hat{A}$, so
$\hat{m}(\omega\otimes\eta )=\omega\eta$, for $\omega,\eta\in\hat{A}$, which means that
\[\hat{m}(\hat{a}\otimes\hat{b})(c) =(\hat{a}\hat{b})(c)=(\hat{a}\otimes\hat{b})\D (c)
=(\varphi\otimes\varphi )(\D (c)(a\otimes b)),\]
for $a,b,c\in A$, and remains valid also for $c\in M(A)$.

\bdefin \label{d-monoid}
A semigroup in a (strict) tensor category $\2C$ is a pair $(Q,m)$, where $Q$ is an object and 
$m: Q\otimes Q\rarr Q$ satisfies $m\mcirc(m\otimes\id_Q)=m\mcirc(\id_Q\otimes m)$. A monoid is a
triple $(Q,m,\eta)$ where $(Q,m)$ is a semigroup and $\eta: \11\rarr Q$ satisfies
$m\mcirc(\eta\otimes\id_Q)=m\mcirc(\id_Q\otimes\eta)=\id_Q$. 
Two semigroups (monoids) $(Q,m,\eta), \ (Q',m',\eta')$ are isomorphic if there exists an
isomorphism $s: Q\rarr Q'$ such that $s\circ m=m'\circ(s\otimes s)$ (and $s\circ\eta=\eta'$).
\edefin

\bprop
\label{p-formon}
Let notation be as above and consider 
the linear map $\tilde{m} =\2F^{-1}\hat{m}(\2F\otimes\2F ):A\otimes A\rarr A$. Then:
\begin{enumerate}
\item $\varphi (c\tilde{m}(x))=(\varphi\otimes\varphi )(\D (c)x)$ for $x\in A\otimes A$ and
$c\in M(A)$. 
\item $\tilde{m}(a\otimes b)=(\varphi\otimes\iota )[((S^{-1}\otimes\iota )\D (b))(a\otimes 1)]
=(\iota\otimes\varphi )[(1\otimes S^{-1}(b))\D (a)]$ for $a,b\in A$.
\item $\tilde{m}(\tilde{m}\otimes\iota )=\tilde{m}(\iota\otimes \tilde{m} )$, so $\tilde{m}$ is a 
multiplication on $A$.
\item $\tilde{m}(\D(a)x)=a\tilde{m}(x)$ for $a\in A$ and $x\in A\otimes A$.
\end{enumerate}
\eprop

\prf
The identity $\2F\tilde{m}=\hat{m}(\2F\otimes\2F )$ means that 
\[\varphi (c\tilde{m}(a\otimes b))
=(\2F\tilde{m}(a\otimes b))(c)
=(\hat{m}(\2F\otimes\2F )(a\otimes b))(c)
=\hat{m}(\hat{a}\otimes \hat{b})(c)
=(\varphi\otimes\varphi )(\D (c)(a\otimes b)),\]
for $a,b,c\in A$, which proves statement 1.
Statement 2 is now immediate from Lemma \ref{l-strinv} and faithfulness of $\varphi$.
To show 3., calculate
\begin{eqnarray*}   \tilde{m}(\tilde{m}\otimes\iota ) &=&
\2F^{-1}\hat{m}(\2F\otimes\2F)(\2F^{-1}\hat{m}(\2F\otimes\2F)\otimes \iota ) 
  = \2F^{-1}\hat{m}(\hat{m}\otimes\iota)(\2F\otimes\2F\otimes\2F )  \\
 &= & \2F^{-1}\hat{m}(\iota\otimes\hat{m})(\2F\otimes\2F\otimes\2F ) =\tilde{m}(\iota\otimes \tilde{m} ).
\end{eqnarray*}
Claim 4 is checked by using 1. and computing
\[ \varphi (c\tilde{m}(\D (a)x))=(\varphi\otimes\varphi )(\D (c)\D (a)x)
  =(\varphi\otimes\varphi )(\D (ca)x)=\varphi (ca\tilde{m}(x)), \] 
for $x\in A\otimes A$ and $a,c\in A$. Now 4. follows by faithfulness of $\varphi$. 
\qed

\bcoro We have $\tilde{m}\in\Hom(\pi_l\times\pi_l,\pi_l)$,  
and $(\pi_l,\tilde{m})$ is a semigroup in $\Rep(A,\D)$.
\ecoro

\prf By the previous proposition the linear map 
$\tilde{m}:A\otimes A\rarr A$ is associative and satisfies
\[ \tilde{m}(\pi_l\times\pi_l )(a)x=\tilde{m}(\D(a)x)=\pi_l (a)\tilde{m}x, \] 
for $a\in A$ and $x\in A\otimes A$. Thus $\tilde{m}$ is an intertwiner from
$\pi_l\times\pi_l$ to $\pi_l$. 
\qed 
 
\brem If $s:\pi\rarr\pi'$ is bounded w.r.t.\ the scalar products on $K, K'$, 
then $s^*$ can be defined as
the adjoint of the unique extension of $s$ to the Hilbert space completions. Therefore, the reader
might wonder why we do not work with the usual tensor $*$-category of non-degenerate
$*$-representations of a discrete AQG on Hilbert spaces. Considering bounded morphisms is,
however, not sufficient for our purposes, since the morphism $\tilde{m}$, which plays a fundamental
role in our considerations, is not bounded w.r.t the 2-norms on $H\otimes H$ and $H$.
To see this it suffices to consider the simple case of 
$(\hat{A},\hat{\D})$ with $\hat{A}_r =C (\7T )$, so $\hat{m}(f\otimes g)(s,t)=f(s)g(t)$, for 
$f,g\in C(\7T )$ and $s,t\in\7T$. Here $(\hat{A}_r ,\hat{\D}_r )$ is the analytic extension of 
$(\hat{A},\hat{\D})$ in the sense of \cite{KD}, so $\hat{A}_r$ is a unital C*-algebra and
$(\hat{A}_r ,\hat{\D}_r )$ is a compact quantum group in the sense of S.L. Woronowicz. By $C (\7T )$
we then mean the unital C*-algebra of all continuous complex valued 
functions of the circle $\7T$ with pointwise algebraic operations and uniform norm.
Since $\2F$ is an isometry by Plancherel's formula, we must then require $\hat{m}$ to be bounded 
w.r.t. the 2-norms on the Hilbert spaces $L^2 (\7T )$ and $L^2 (\7T\times\7T )$ of square integrable
functions on $\7T$ and $\7T\times\7T$ (obtained from the GNS-constructions of $\hat{A}$ and
$\hat{A}\otimes\hat{A}$ described in \cite{KD}), and this is clearly false.   
Thus one cannot define $\tilde{m}^*$ by extension to the Hilbert space completion. Also Proposition
\ref{p-rep*} is not applicable, since in general $\pi_l\times\pi_l$  
is not in $\Rep_*(A,\D)$. 
\erem

\bprop
\label{p-monoid}
Let $(A,\D )$ be an AQG. Then $\Hom(\ve,\pi_l)\ne \{0\}$ iff $(A,\D )$ is discrete.
In this case, the map $\tilde{\eta}: c\mapsto c\2F^{-1}(1_{\hat{A}})$ belongs to 
$\Hom(\ve,\pi_l)$ and  $(\pi_l,\tilde{m},\tilde{\eta})$ is a monoid, which we call the regular
monoid. We have $\tilde{\eta}c=cI_0$. Since $\tilde{\eta}: \ve\rarr\pi_l$ is a morphism in
$\Rep_*(A,\D)$,  the adjoint $\tilde{\eta}^*$ exists and $\tilde{\eta}^*=\ve$. 
\eprop
\prf
For every morphism $\eta\in\Hom(\ve,\pi_l)$ we have 
\[ \varepsilon (a)\eta(1)=\eta (\varepsilon (a)1)=\pi_l (a)\eta (1)=a\eta (1) \]
for $a\in A$, saying that $\eta(1)$ is a left integral in $A$. Thus $\eta\to \eta (1)$ is a bijection
from $\Hom(\ve,\pi_l)$ to the space of left integrals in $A$. By definition an AQG $(A,\D)$ 
is discrete iff a non-zero left integral exists, and in this case it is unique up to a scalar. 

If $(A,\D)$ is discrete then $(\hat{A},\hat{\D})$ is compact, i.e.\ $\hat{A}$ has a 
unit $1_{\hat{A}}$. Now
\[ \2F(I_0)(a)=\varphi(aI_0)=\ve(a)\varphi(I_0)=\ve(a). \]
Thus $\2F(I_0)=\varepsilon =1_{\hat{A}}$ and $\tilde{\eta}(1)=I_0$, which is a left integral in $A$,
so $\tilde{\eta}\in\Hom(\ve,\pi_l)$. 

Finally, the equalities
\[ (\tilde{\eta}(c),a)_A=\varphi(a^*cI_0)=c\ve(a^*)=c\ol{\ve(a)}=(c,\ve(a))_\7C , \]
for $c\in\7C$ and $a\in A$, show that $\tilde{\eta}^*=\ve$.
\qed

\brem 1. The above result shows in particular that a monoid structure on the regular representation
exists only if $(A,\D)$ is discrete. It turns out that the multiplication $\tilde{m}$ is in 
general not unique, not even up to isomorphisms of $\pi_l$.

2. If $(A,\D)$ is a discrete and quasitriangular AQG with $R$-matrix $R$, the categories
$\Rep_f(A,\D)$ and $\Rep(A,\D)$ are braided. It is therefore natural to ask whether the monoid
$(\pi_l,\tilde{m},\tilde{\eta})$ is commutative in the sense that 
$\tilde{m}\circ c_{\pi_l,\pi_l}=\tilde{m}$, where $c$ denotes the braiding.
One can easily show that this is the case iff 
$R=1\otimes 1$. In that case, $(A,\D)$ is cocommutative and the representation
categories are symmetric.
\erem

\bdefin
A comonoid in a (strict) tensor category $\2C$ is a triple $(Q,\D,\varepsilon )$, where $Q$ is an
object and $\D: Q\rarr Q\otimes Q,\ \varepsilon : Q\rarr \11$ satisfy
$\D\otimes\id_Q \mcirc \D\ =\id_Q\otimes \D \mcirc \D$ and
$\varepsilon\otimes\id_Q\mcirc \D=\id_Q\otimes\varepsilon \mcirc \D=\id_Q$. 
\edefin

For a compact AQG we have the following easy result.

\bprop
\label{t-compcomon} Let $(A,\D)$ be an AQG. The map $\ve: A\rarr\7C$ is in $\Hom(\pi_l,\ve)$. 
Furthermore, $(A,\D)$ is compact iff $\D(A)\subset A\otimes A$ iff 
$\D\in\Hom(\pi_l,\pi_l\times\pi_l)$. In this case $(\pi_l, \D,\ve)$ is a comonoid in $\Rep(A,\D)$,
which we call the regular comonoid. 
\eprop

\prf
For any AQG we have the equation $m(S\otimes\iota)\Delta(a)=\ve(a)I$ in $M(A)$. If 
$\D(A)\subset A\otimes A$ the left hand side and therefore the unit $I$ belongs to $A$. The
remaining facts are obvious consequences of $\ve$ and $\D$ being algebra homomorphisms.
\qed

\brem
Again, one might try to work with the usual tensor $*$-category of unital $*$-representations of a
compact AQG $(A,\D )$ on Hilbert spaces. There is no problem with $\D$, as it is an isometry, but
$\varepsilon :A\rarr\7C$ has in the general case no continuous extension w.r.t the 2-norm on $A$
given by the GNS-construction. 
\erem

\bdefin
A Frobenius algebra in a tensor category $\2C$ is a quintuple $(Q,m,\eta,\D,\varepsilon)$
such that $(Q,m,\eta )$ is a monoid in $\2C$, $(Q,\D,\varepsilon )$ is a comonoid in $\2C$, and the
following compatibility condition holds
\begin{equation} \label{e-frob}
\id_Q\otimes m\mcirc\D\otimes \id_Q =\D\mcirc m =m\otimes\id_Q \mcirc\id_Q\otimes\D.
\end{equation}
\edefin

\bprop \label{p-frob}
Let $(A,\D )$ be a finite dimensional AQG. Then $\tilde{m}^*=\D$ so the regular monoid
and comonoid are each others adjoints:
$(\pi_l,\tilde{m},\tilde{\eta})^*\equiv(\pi_l,\tilde{m}^*,\tilde{\eta}^*)=(\pi_l,\D,\ve)$. 
Furthermore, the quintuple $(\pi_l,\tilde{m},\tilde{\eta},\D,\ve)$ is a Frobenius algebra in
$\Rep_f(A,\D)$, which we call the regular Frobenius algebra.
\eprop

\prf That $\tilde{\eta}^*=\ve$ is shown in Proposition \ref{p-monoid}. If $(A,\D)$ is finite
dimensional, $\tilde{m}^*=\D$ follows from statement 1 of Proposition \ref{p-formon}. The Frobenius
property (\ref{e-frob}) will be shown at the end of the next subsection. (Cf.\ also \cite{mue09}.)
\qed

\brem 1. Conversely, the existence of both the regular monoid and the regular comono\-id requires
$(A,\D)$ to be discrete and compact, thus $A$ is finite dimensional.

2. Given a Frobenius algebra $(Q,m,\eta ,\D,\varepsilon )$, it is easy to show that the morphisms  
$\varepsilon\circ m : Q\otimes Q\rarr\11$ and $\D\circ\eta :\11\rarr Q\otimes Q$ satisfy the 
triangular equations \cite{kas}, i.e.\ the object $Q$ is its own two-sided dual. If $\2C$ is a
$*$-category and $m^*=\D,\ \eta^*=\ve$ we have $(\D\circ\eta)^*=\ve\circ m$ and we obtain a solution
of the conjugate equations \cite{LR}. In the case considered above, this in particular implies that
$\pi_l$ is a finite dimensional object in $\Rep(A,\D)$, thus again $A$ is finite dimensional.
\erem

We show now how one can recover the intrinsic group, cf.\ \cite{MRT}, from the regular monoid of a
discrete AQG.

\bdefin
Let $\2C$ be a tensor category and $\2C_*$ a full $*$-subcategory. Let $(Q,m, \eta )$ be a monoid in
$\2C$ with $Q\in\2C_*$. Denote by $G_Q$ the group in $\End Q$ given by
\[G_Q =\{ t\in\End\,Q\ |\ t\circ t^* =t^*\circ t=\id_Q ,\ \ m\mcirc t\otimes t=t\circ m\} \]
with group multiplication being composition of arrows,
so the unit of $G_Q$ is $\id_Q$ and the inverse $t^{-1}$ of $t\in G_Q$ is $t^*$. 
The group $G_Q$ is called the intrinsic group of the monoid $(Q,m, \eta )$.
\edefin

\bprop
Let $(A,\D )$ be a discrete AQG with intrinsic group $G$ defined by
\[ G=\{g\in M(A)\ |\ \D g=g\otimes g,\ g^* g =gg^* =I\}, \]   
which is compact w.r.t.\ the product topology on $M(A)$.
Let $G_{\pi_l}\subset \Hom(\pi_l,\pi_l)$ be the intrinsic group of 
the regular monoid $(\pi_l, \tilde{m}, \tilde{\eta})$ with 
topology defined in Proposition \ref{p-rep*}. 
Then $G\cong G_{\pi_l}$ as topological groups. 
\eprop

\prf
By Proposition \ref{p-disc} $M(A)\cong R(A)$ as unital algebras, and $R(A)$ and 
$\Hom (\pi_l ,\pi_l )$ are anti-isomorphic as unital $*$-algebras. 
Let $\rho :A\rightarrow A$ denote the linear map such that
$\varphi (ab)=\varphi (b\rho (a))$ for $a,b\in A$.  
It suffices to show that for any $g\in R(A)$, we have 
$\tilde{m}\circ (g\otimes g)=g\circ\tilde{m}$ iff $\D g=g\otimes g$.
But $\tilde{m}\circ (g\otimes g)=g\circ\tilde{m}$ means that
$\tilde{m}\circ (g\otimes g)(a\otimes b)=g\circ\tilde{m}(a\otimes b)$, for $a,b\in A$,
or $\tilde{m}(ag\otimes bg)=\tilde{m}(a\otimes b)g$,  
which by faithfulness of $\varphi$, can be expressed as
\[\varphi (c\tilde{m}(ag\otimes bg))=\varphi (c\tilde{m}(a\otimes b)g)
=\varphi (\rho^{-1}(g)c\tilde{m}(a\otimes b)),\]
for $a,b,c\in A$.
Hence by Proposition \ref{p-formon} the formula $\tilde{m}\circ (g\otimes g)=g\circ\tilde{m}$ 
can be rewritten as
\begin{eqnarray*} \lefteqn{ 
   (\varphi\otimes\varphi )(\D c(ag\otimes bg))=(\varphi\otimes\varphi )(\D ({\rho}^{-1}(g)c)(a\otimes b)) } \\
 && =(\varphi\otimes\varphi )(\D \rho^{-1}(g)\D c(a\otimes b))
=(\varphi\otimes\varphi )(\D (c)(a\otimes b)(\rho\otimes\rho )\D \rho^{-1}(g)),
\end{eqnarray*}
for $a,b,c\in A$. Thus again by faithfulness of $\varphi$, we see that
$\tilde{m}\circ (g\otimes g)=g\circ\tilde{m}$ iff $(\rho\otimes\rho )\D\rho^{-1}(g)=g\otimes g$.

We assert now that $(\rho\otimes\rho )\D\rho^{-1} =\D$ for any discrete AQG, 
which clearly completes the proof of the proposition. Since 
$(S^2\otimes \rho )\D=\D\rho$ for any AQG, we see that 
$(\rho\otimes\rho )\D\rho^{-1} =\D$ iff $\rho=S^2$, which holds for discrete AQG.
\qed


\subsection{The Absorbing Property}
The following obvious fact will be used without further reference.

\blemma Let $(A,\D)$ be a AQG. If $\theta$ is a $*$-representation of $(A,\D)$ on $K$ then 
\[\begin{array}{ccc}
  (\theta\times\pi_l)(a)(v\otimes x)=\sum_i \theta(a_i)v\otimes b_i &\quad\mbox{where}
    \quad &  \sum_i a_i\otimes b_i = \D(a)(1\otimes x), \\
  (\pi_l\times\theta)(a)(x\otimes v) = \sum_i a_i\otimes\theta(b_i)v &\quad\mbox{where}\quad    
   &\sum_i a_i\otimes b_i = \D(a)(x\otimes 1)
\end{array}\]
for $v\in K$ and $x\in A$.
\elemma

\bprop 
\label{p-abs}
Let $(A,\D)$ be an AQG. For every $\theta\in\Rep(A,\D)$ we have the absorption property 
$$\DS \theta\times\pi_l\ \cong\ \pi_l\times\theta\ \cong\  I_{\theta}\times\pi_l\ \cong\ \dim K\,\pi_l$$
for $\pi_l$, where $I_{\theta}$ is the $*$-representation of $A$ on $K$ given by 
$I_{\theta}(a)=\varepsilon (a)\id_K$ for $a\in A$.  
\eprop
\prf We start by showing
$\DS \theta\times\pi_l\cong I_{\theta}\times\pi_l$.
Define a linear map 
$U_{\theta} :K\otimes A\rarr K\otimes A$ by
\[U_{\theta}(\theta (a)v\otimes x)=\sum_i \theta (a_i )v\otimes x_i ,\]
where $\sum_i a_i\otimes x_i =\D (x)(a\otimes 1)$ for $a,x\in A$. To see that $U_{\theta}$
is well-defined suppose $\sum_j\theta (a^j )v^j\otimes x^j =0$, where $a^j ,x^j\in A$ and $v^j\in
K$, and write $\sum_j \D (x^j )(a^j\otimes 1)=\sum_{ij} a^j_i\otimes x^j_i$ with $a^j_i ,x^j_i\in A$. 
We must show that $\sum_{ij} \theta (a^j_i )v^j\otimes x^j_i =0$, and in doing so,
we can assume that $(x^j )$ are linear independent, so $\theta (a^j )v^j =0$ for all $j$.
Pick a two sided local unit $e$ for the collection $\{a^j_i ,x^j _i\}$ and a local unit $e'$ for $e$.
Then 
\begin{eqnarray*}    \sum_{ij} \theta (a^j_i )v^j\otimes x^j_i &=& 
   \sum_{ij} (\theta (e)\otimes e)(\theta (a^j_i) v^j\otimes x^j_i e')(v^j\otimes e) \\
 &=& (\theta (e)\otimes e)\sum_j (\theta\otimes\iota )(\D (x^j )(a^j\otimes 1))(v^j\otimes e'e) \\
 &=& (\theta (e)\otimes e)\sum_j (\theta\otimes\iota )(\D (x^j )(1\otimes e'))(\theta (a^j)v^j\otimes e)
   =0
\end{eqnarray*}
as $\theta (a^j )v^j =0$ for all $j$. Thus $U_\theta$ is well-defined.

Furthermore, for $a,b,x\in A$ and $v\in K$, we have  
\[(\theta\times \pi_l )(b)U_{\theta} (\theta (a)v\otimes x)
=\sum_i (\theta\times \pi_l )(b)(\theta (a_i )v\otimes x_i )
=\sum_{ik} \theta (b_k^i )\theta (a_i )v\otimes y_k^i =\sum_{ik} \theta (b_k^i a_i )v\otimes y_k^i ,
\] 
where $\sum_i a_i\otimes x_i =\D (x)(a\otimes 1 )$ and 
$\sum_k b_k^i\otimes y_k^i =\D (b)(1\otimes x_i )$. But 
\[\sum_{ik} b_k^i a_i\otimes y_k^i =\sum_i \D (b)(1\otimes x_i )(a_i\otimes 1)
  =\D(b)\D(x)(a\otimes 1) =\D (bx)(a\otimes 1),\]
so 
$(\theta\times \pi_l )(b)U_{\theta} (\theta (a)v\otimes x)
=U_{\theta}(\theta (a)v\otimes bx)$
for $a,b,x\in A$ and $v\in K$.
On the other hand, if we write $\D (b)(1\otimes x)=\sum_i c_i\otimes z_i$, for $b,x\in A$, and
calculate 
\begin{eqnarray*} \lefteqn{
(I_{\theta}\times \pi_l )(b)(\theta (a)v\otimes x)
= \sum_i I_{\theta}(c_i )\theta (a)v\otimes z_i
=\sum_i \varepsilon (c_i )\theta (a)v\otimes z_i } \\
  &&= \theta (a)v\otimes (\sum_i \varepsilon (c_i )z_i )
=\theta (a)v\otimes (\varepsilon\otimes\iota )[\D (b)(1\otimes x)]
=\theta (a)v\otimes bx ,
\end{eqnarray*}
for $a\in A$ and $v\in K$, we see that
\[(\theta\times \pi_l )(b)U_{\theta} (\theta (a)v\otimes x)
=U_{\theta}(\theta (a)v\otimes bx)=U_{\theta}(I_{\theta}\times \pi_l )(b)(\theta (a)v\otimes x),\]
for $a,b,x\in A$ and $v\in K$, so
\[(\theta\times \pi_l )(b)U_{\theta}=U_{\theta}(I_{\theta}\times \pi_l )(b),\]
for $b\in A$ and therefore $U_{\theta}\in\Hom (I_{\theta}\times \pi_l ,\theta\times \pi_l )$
in $\Rep(A,\D )$.

By the cancellation laws for $(A,\D )$, we see that $U_{\theta}$ is invertible, and thus 
$\DS \theta\times\pi_l\cong I_{\theta}\times\pi_l$.
In fact, the inverse of $U_{\theta}$ in $\Rep(A,\D)$ is given by the following formula
$U_{\theta}^{-1} (\theta (a)v\otimes x) =\sum_i \theta (a_i )v\otimes x_i$,
where
\[ \sum_i a_i\otimes x_i = ((S^{-1}\otimes\iota)\D (x))(a\otimes 1) \]
for $a,x\in A$ and $v\in K$.

Similarly, one shows that the linear map $V_{\theta} :A\otimes K\rarr A\otimes K$ given by
\[ V_{\theta}(x\otimes \theta (a)v)=\sum_i x_i\otimes\theta (a_i )v, \]
for $a,x\in A$ and $v\in K$, and where $\sum_i x_i\otimes a_i =\D (x)(1\otimes a)$,
is well-defined and is an isomorphism in $\Hom (\pi_l\times I_{\theta},\pi_l\times\theta )$. 
Thus $\pi_l\times I_{\theta}\cong\pi_l\times\theta$ in $\Rep(A,\D)$.

From the previously derived expression
\[(I_{\theta}\times\pi_l )(b)(\theta (a)v\otimes x)=\theta (a)v\otimes bx,\]
valid for $a,b,x\in A$ and $v\in K$, 
we see that 
$(I_{\theta}\times\pi_l )(b) =I_{B(K)}\otimes\pi_l (b)$, so 
$I_{\theta}\times\pi_l\cong (\dim K)\pi_l$ and similarly $\pi_l\times I_{\theta}\cong (\dim K)\pi_l$. 
\qed

\brem 
In the case of a discrete AQG, where $\Rep(A,\D)$ is semisimple,
Proposition \ref{p-abs} can also be proven using 1. in Proposition \ref{p-disc}
together with Proposition \ref{l-absor} below.
\erem

\bprop Let $(A,\D)$ be a discrete AQG and $\theta\in\Rep_f(A,\D)$. Then 
$\pi_l\times\theta, \pi_l\times I_\theta,\theta\times\pi_l, I_\theta\times\pi_l\in\Rep_*(A,\D)$, and
the morphisms $U_\theta,\ V_\theta$ considered in the preceding proposition are unitary.
\eprop

\prf In the discrete case $\pi_l\in\Rep_*(A,\D)$ by Proposition \ref{p-disc} 
and the same is true for $\pi_l\times\theta$ and
$\theta\times\pi_l$. Since we know that $U_\theta$ is invertible, we need only show that it is an 
isometry. (As always, the regular representation $\pi_l$ is understood to come with its scalar
product defined using $\varphi$.) Recall that $U_\theta$ is defined by
\[U_{\theta}(\theta (a)v\otimes x)=\sum_i \theta (a_i )v\otimes x_i , \]
where $\sum_i a_i\otimes x_i =\D (x)(a\otimes 1)$ for $a,x\in A$. We compute
\begin{eqnarray*} \lefteqn{ 
(U_{\theta}(\theta (a)v\otimes x),U_{\theta}(\theta (a)v\otimes x))= 
  \sum_{ij}(\theta (a_i )v\otimes x_i ,\theta (a_j )v\otimes x_j ) } \\
&&  =\sum_{ij}(\theta (a_i )v,\theta (a_j )v )\varphi (x_j^* x_j ) 
  =(\theta (\sum_{ij} a_j ^*a_i \varphi (x_j^* x_j ))v,v ) \\
&&  =(\theta (\sum_{ij} (\iota\otimes\varphi )(a_j ^*a_i\otimes x_j^* x_j ))v,v )
  =(\theta ((\iota\otimes\varphi )([\D (x)(a\otimes 1)]^* [\D (x)(a\otimes 1)]))v,v )  \\
&&    =(\theta ((\iota\otimes\varphi )((a^*\otimes 1)\D (x^* x)(a\otimes 1)))v,v )
   =(\theta (a^* a\varphi(x^* x))v,v ) \\
&&  =(\theta (a)v ,\theta (a)v)(x ,x ) = (\theta (a)v\otimes x,\theta (a)v\otimes x ),
\end{eqnarray*}
for $a,x\in A$ and $v\in K$. Thus $U_{\theta}$ is an isometry. The same is true for $V_\theta$
by a similar computation.
\qed

\brem 1. Let $\theta$ be a $*$-representation of $(A,\D)$ with $U_\theta$ as above. Clearly, there
is a $*$-representation $\tilde{\theta}$ of $(A,\D_\op)$ which coincides with $\theta$ as a map from
$A$ to $\End\,K$. It is then easy to see that $V_{\tilde{\theta}}=\Sigma U_\theta\Sigma^{-1}$, where 
$\Sigma: K\otimes A\rarr A\otimes K$ is the flip map. This observation obviates separate proofs for
$V_\theta$. 

2. Note that the assumptions on $(A,\D)$ and $\theta$ were only made in order for $U_\theta^*$ to be
definable in $\Rep_*(A,\D)$. The computation showing that $U_\theta$ is isometric holds in general
and provides an alternative proof for the well-definedness of $U_\theta$.
\erem

\bprop 
\label{p-functorial}
The morphisms $V_\theta: \pi_l\times I_\theta\rarr \pi_l\times\theta$ are natural w.r.t.\
$\theta$, i.e.\ the diagrams
\[ \begin{diagram} \pi_l\times I_\theta & \rTo^{V_\theta} & \pi_l\times \theta \\
\dTo^{\id_{\pi_l}\otimes s} & & \dTo_{\id_{\pi_l}\otimes s} \\
\pi_l\times I_{\theta'} & \rTo^{V_{\theta'}} & \pi_l\times \theta'
\end{diagram}\]
commute for all $s:\theta\rarr\theta'$, and similarly for $U_\theta$.
\eprop

\prf
This is obvious by definition of $V_\theta$.
\qed

Having defined monoids in tensor categories, we will also need the notion of a module over a
monoid.

\bdefin \label{d-module}
Let $\2C$ be a tensor category and $(Q,m)$ a semigroup in $\2C$. Then a (left) $Q$-module is a pair 
$(X,\mu)$, where $X\in\2C$ and $\mu: Q\otimes X\rarr X$ satisfies 
\[ \mu\mcirc m\otimes\id_X=\mu\mcirc\id_Q\otimes\mu. \]
For a monoid $(Q,m,\eta)$ we require in addition that $\mu\mcirc\eta\otimes\id_X=\id_X$.
With
\[ \Hom_{Q\!-\!\Mod}((X,\mu),(X',\mu'))=\{ s\in\Hom(X,X')\ | \ s\mcirc\mu=\mu'\mcirc\id_Q\otimes s\}
\]
as morphisms, the $Q$-modules form a category which we denote by $Q\!-\!\Mod$.
\edefin

\bprop \label{p-Vmod}
The diagram
\[ \begin{diagram} \pi_l\times \pi_l\times I_\theta & \rTo^{\tilde{m}\otimes\id_\theta} & \pi_l\times I_\theta \\
  \dTo^{\id_{\pi_l}\otimes V_\theta} && \dTo_{V_\theta} \\
  \pi_l\times \pi_l\times \theta & \rTo^{\tilde{m}\otimes\id_\theta} & \pi_l\times \theta
\end{diagram} \]
commutes. We have similar commutative diagrams for the morphisms
$V_\theta^*: \pi_l\times \theta\rarr\pi_l\times I_\theta$. 
\eprop

\prf We must show that 
$(\tilde{m}\otimes\iota )(\iota\otimes V_\theta)=V_\theta(\tilde{m}\otimes\iota)$ 
as maps from $A\otimes A\otimes K$ to $A\otimes K$. Let $a,c,x,y\in A$ and $v\in K$ and write 
$\D (x)(1\otimes a)=\sum_i x_i\otimes a_i$, where $a_i ,x_i\in A$. Then 
\bean \lefteqn{ 
  (\varphi c\otimes\iota )(\tilde{m}\otimes\iota )(\iota\otimes V_\theta )(y\otimes x\otimes \theta (a)v)
 = \sum_i \varphi (c\tilde{m}(y\otimes x_i ))\theta (a_i )v } \\
  &&= \theta \left( \sum_i (\varphi\otimes\varphi )(\D (c)(y\otimes x_i ))a_i \right)v \\
 &&= \theta ((\varphi\otimes\varphi\otimes\iota )((\D (c)\otimes 1)(1\otimes \D (x))(y\otimes 1\otimes a)))v ,
\eean
whereas if we write
\[ \D (\varphi\otimes\iota )[((S^{-1}\otimes\iota )\D (x))(y\otimes 1)](1\otimes a)=\sum_j y_j\otimes b_j\]
for $b_j ,y_j\in A$ and use 2. in Proposition \ref{p-formon}, we get
\begin{eqnarray*}
  (\varphi c\otimes\iota )V_\theta (\tilde{m}\otimes\iota )(y\otimes x\otimes \theta (a)v)
  &=& 
(\varphi c\otimes\iota )V_\theta ((\varphi\otimes\iota )[((S^{-1}\otimes\iota )\D (x))(y\otimes
1)]\otimes \theta (a)v) \\
  &=&\theta \left( \sum_j \varphi (cy_j )b_j \right) v. 
\end{eqnarray*}
Hence 
$(\tilde{m}\otimes\iota )(\iota\otimes V_\theta )=V_\theta (\tilde{m}\otimes\iota )$
follows if
\[(\varphi\otimes\varphi\otimes\iota )((\D (c)\otimes 1)(1\otimes \D (x))(y\otimes 1\otimes a)=\sum_j \varphi (cy_j )b_j .\] 
Now,
\[\sum_j \varphi (cy_j )b_j =(\varphi c\otimes\iota )(\sum_j y_j\otimes b_j )
=(\varphi c\otimes\iota )[\D (\varphi\otimes\iota )[((S^{-1}\otimes\iota )\D (x))(y\otimes 1)](1\otimes a)] ,\]
so $(\tilde{m}\otimes\iota )(\iota\otimes V_\theta )=V_\theta (\tilde{m}\otimes\iota )$ if
\[(\varphi c\otimes\iota )\D (\varphi\otimes\iota )[((S^{-1}\otimes\iota )\D (x))(y\otimes 1)]
=(\varphi\otimes\varphi\otimes\iota )((\D(c)\otimes 1)(1\otimes\D (x))(y\otimes 1\otimes 1)).\]
But
\bean L.H.S. &=& (\varphi c\otimes\iota )\D (\hat{y}\otimes\iota )(S^{-1}\otimes\iota )\D (x)
=(\varphi c\otimes\iota )\D (\hat{y}S^{-1}\otimes\iota )\D (x) \\
 &=& (\hat{y}S^{-1}\otimes\varphi c\otimes\iota )(\iota\otimes\D )\D (x)
=(\hat{y}S^{-1}\otimes\varphi c\otimes\iota )(\D\otimes\iota )\D (x) ,\eean
whereas by strong left-invariance of $\varphi$, we get
\bean R.H.S. &=& (\hat{y}\otimes\varphi\otimes\iota )((\D (c)\otimes 1)(1\otimes\D (x)))
=(\hat{y}S^{-1}S\otimes\varphi\otimes\iota )((\D (c)\otimes 1)(1\otimes\D (x))) \\
 &=& (\varphi\otimes\iota )(c\otimes 1((\hat{y}S^{-1}\otimes\iota )\D\otimes\iota )\D (x))
=(\hat{y}S^{-1}\otimes\varphi c\otimes\iota )(\D\otimes\iota )\D (x),\eean
as desired. Now, replacing $V_\theta$ by its inverse $V_\theta^*$, the direction of the vertical
arrows in the diagram is reversed, and we see that also $V_\theta^*$ is a $\pi_l$-module morphism. 
\qed

\bcoro
The morphisms $V_\theta: \pi_l\times I_\theta\rarr \pi_l\times \theta$ and 
$V_\theta^*: \pi_l\times \theta\rarr\pi_l\times I_\theta$ are morphisms of $\pi_l$-modules.
\ecoro

\bprop \label{p-qf}
Let $(A,\D)$ be a discrete AQG. Then 
\begin{equation} \label{e-qf}
 (\tilde{m}\otimes\iota)(y\otimes(\D(x)(1\otimes a)))=\D(\tilde{m}(y\otimes x))(1\otimes a) 
\end{equation}
holds for $x,y,a\in A$.
\eprop

\prf First note that for $V_\theta$ with $\theta =\pi_l$, we
have $V_\theta (x\otimes a)=\D (x)(1\otimes a)$ for $a,x\in A$. To see this write 
$\D (x)(1\otimes a)=\sum_i x_i\otimes a_i$ with $a_i ,x_i\in A$ and pick a right-sided local unit
$e\in A$ for $\{a, a_i\}$. Then 
\[V_\theta (x\otimes a)=V_\theta (x\otimes \pi_l (a)e)=\sum_i x_i\otimes\pi_l (a_i )e =\sum_i x_i\otimes a_i e 
  =\sum_i x_i\otimes a_i =\D (x)(1\otimes a).\]
By Proposition \ref{p-Vmod} we have 
$(\tilde{m}\otimes\iota )(\iota\otimes V_\theta )=V_\theta (\tilde{m}\otimes\iota )$. Thus
\bean \lefteqn{  (\tilde{m}\otimes\iota)(y\otimes(\D(x)(1\otimes a)))
=(\tilde{m}\otimes\iota )(\iota\otimes V_\theta )(y\otimes x\otimes a) } \\
 && = V_\theta (\tilde{m}\otimes\iota )(y\otimes x\otimes a)
=\D(\tilde{m}(y\otimes x))(1\otimes a) \eean
for $x,y,a\in A$. 
\qed

\brem
Clearly, it follows from the proof of this proposition that $\tilde{m}$ has the property stated in the
proposition iff $V_{\pi_l}$ is a $\pi_l$-module map. In the discrete case, where every
representation is a direct sum of representations contained in $\pi_l$, Proposition \ref{p-Vmod}
can therefore also be deduced using the naturality property of Proposition \ref{p-functorial}.
\erem

\noindent{\it End of proof of Proposition \ref{p-frob}.}
If $(A,\D)$ is finite dimensional, we can put $a=1$ in (\ref{e-qf}) and obtain
\[ (\tilde{m}\otimes\iota)(\iota\otimes\D)(y\otimes x)=\D(\tilde{m}(y\otimes x)) \]
for $x,y\in A$. In categorical terms this is the equality 
$\tilde{m}\otimes\id_{\pi_l}\,\circ\,\id_{\pi_l}\otimes\D=\D\circ\tilde{m}$ in
$\End(\pi_l\otimes\pi_l)$. Using the $*$-operation and $\Delta^*=\tilde{m}$ we also find  
$\id_{\pi_l}\otimes\tilde{m}\,\circ\,\D\otimes\id_{\pi_l}=\D\circ\tilde{m}$. This completes the
proof of the Frobenius property in Proposition \ref{p-frob}. \qed

We close this section by summarizing the results on the regular representation.

\btheor
Let $(A,\D)$ be an AQG with left regular representation $\pi_l$. Then there exists a morphism 
$\tilde{m}:\pi_l\times\pi_l\rarr\pi_l$ such that $(\pi_l,\tilde{m})$ is a semigroup in the tensor
category $\Rep(A,\D)$. The representation $\pi_l$ has the absorbing property 
$\pi_l\times\theta\ \cong\ \pi_l\times I_{\theta} \cong\ \dim K\,\pi_l$ w.r.t.\
a natural family of equivalences
$V_\theta: \pi_l\times I_\theta\rarr \pi_l\times\theta$ of (left) $\pi_l$-modules. 
Similarly, there are natural equivalences
$U_\theta: I_\theta\times\pi_l\rarr\theta\times\pi_l$ of right $\pi_l$-modules.
These equivalences are unitary whenever $(A,\D)$ is discrete and $\theta\in\Rep_f (A,\D)$.

There exists a morphism $\tilde{\eta}: \ve\rarr\pi_l$ such that $(\pi_l,\tilde{m},\tilde{\eta})$ is
a monoid iff $(A,\D)$ is discrete. In the discrete case, there exists a non-monoidal $*$-subcategory
$\Rep_*(A,\D)\subset\Rep(A,\D)$ containing $\pi_l$.
\etheor


\section{On Monoids, Embedding Functors and AQG}
\subsection{From Monoids to Embedding Functors}
Some of the results in this section will be formulated over any ground field $\7F$.
Let $\mathrm{Vect}_\7F$ denote the tensor category of finite dimensional vector spaces over $\7F$.

\blemma \label{p-faithful}
Let $\2C$ be an $\7F$-linear semisimple category. Then an $\7F$-linear functor 
$F: \2C\rarr\mathrm{Vect}_\7F$ is faithful (i.e.\ $F(s)=0$ for $s: X\rarr Y$ implies $s=0$)
if $F(X)$ is non-zero for every irreducible $X\in\2C$.
\elemma
\prf Suppose $F(X)$ is non-zero for every irreducible $X$ and 
consider $s: X\rarr Y$ such that $F(s)=0$. Let $I_\2C$ be the set of isomorphism 
classes of irreducible objects with chosen representatives $X_i, i\in I_\2C$. 
Let $(v_{i\alpha})$ be bases in $\Hom(X_i,X)$ with dual bases $(v'_{i\alpha})$
 satisfying $v_{i\alpha}'\,\circ\, v_{j\beta}=\delta_{ij}\delta_{\alpha,\beta}\id_{X_i}$ and 
$\sum_\alpha v_{i\alpha}\,\circ\, v_{i\alpha}'=\id_X$. Pick $w_{j\beta}\in\Hom(X_i, Y)$
and $w_{j\beta}'$ similarly.
Since $\Hom(X_i,X_j)=\delta_{ij}\,\id_{X_i}\,\7F$, which implies 
$w_{j\beta}'\,\circ\, s\,\circ\, v_{i\alpha}=\delta_{i,j}\,c_{i\alpha\beta}\,\id_{X_i}$, we can write
\[ s=\sum_{i\alpha, j\beta} w_{j\beta}\,\circ\, w_{j\beta}'\,\circ\, s\,\circ\, v_{i\alpha}\,
\circ\, v_{i\alpha}'
   =\sum_{i\alpha\beta} c_{i\alpha\beta} \ w_{i\beta} \,\circ\, v_{i\alpha}'. \]
Thus 
\[ 0=F(w_{k\eta}')\,\circ\, F(s)\,\circ\, F(v_{k\xi})
= \sum_{i\alpha\beta} c_{i\alpha\beta}\  F(w_{k\eta}'\,\circ\, w_{i\beta}
   \,\circ\, v_{i\alpha}'\,\circ\, v_{k\xi}) 
= c_{k\xi\eta}\,F(\id_{X_k}),\]
for $k,\xi$ and $\eta$. 
By assumption $F(\id_{X_k})\ne 0$ for $k\in I_\2C$, thus all $c_{i\alpha\beta}$ vanish and
$s=0$. 
\qed

\bprop \label{p-embed} 
Let $\2C$ be a semisimple $\7F$-linear tensor category with $\End\11\cong \7F$,
and let $(Q,m,\eta)$ be a monoid in $\hat{\2C}$ such that:  
\begin{enumerate}
\item $\dim\Hom_{\hat{\2C}}(\11,Q)=1$.
\item For every $X\in\2C$, there is an isomorphism $Q\otimes X\cong n(X)Q$ of $Q$-modules
with $n(X)\in\7N$.
\end{enumerate}
Then the functor $E: \ \2C\rarr \Vect_\7F$ defined by
$X\mapsto \Hom_{\hat{\2C}}(\11,Q\otimes X)$ and
\begin{equation} \label{e-Emor}
E(s)\phi=\id_Q\otimes s\mcirc  \phi,  
\end{equation}
where $s:X\rarr Y$ and $\phi\in\Hom(\11,Q\otimes X)$, is a faithful (strong) tensor 
functor with $\dim E(X)=n(X)$.
\eprop

\prf We have $E(X) = \Hom(\11,Q\otimes X)\cong\Hom(\11,n(X)Q)\cong d(X)\Hom(\11,Q)\cong \7F^{n(X)}$,
thus $E(X)$ is a vector space of dimension $n(X)$. Since $E(X)\ne 0$ for every $X\in\2C$, 
Lemma \ref{p-faithful} tells us that $E$ is faithful. 

To see that $E$ is monoidal first observe that $E(\11)=\Hom(\11,Q)=\7F\eta$ by 2. 
Thus there is a canonical isomorphism 
$e: \7F=\11_{\Vect_\7F} \rarr E(\11 )=\Hom(\11,Q)$ defined by $c\mapsto c\eta$. Next we 
define morphisms 
\[ d^E_{X,Y}: E(X)\otimes E(Y)\rarr E(X\otimes Y), \quad \phi\otimes \psi\mapsto
   m\otimes\id_{X\otimes Y}\mcirc\id_Q\otimes\phi\otimes\id_Y\mcirc\psi. \]
In terms of a diagram, this means
\[ d_{X,Y}^E(\phi\otimes\psi)= \quad
\begin{tangle}
\hstep\object{Q}\step[1.5]\object{X}\step\object{Y}\\
\hh\step[-1]\obj{m}\step\hcd\step\id\step\id\\
\hh\id\step\id\step\id\step\id\\
\hh\id\step\frabox{\phi}\step\id\\
\d\step\dd\\
\hh\step\frabox{\psi}
\end{tangle}
\]
By definition (\ref{e-Emor}) of the map $E(s): E(X)\rarr E(Y)$ it is obvious that the family
$(d^E_{X,Y})$ is natural w.r.t.\ both arguments. The equation
\[ d^E_{X_1\otimes X_2,X_3}\mcirc  d^E_{X_1,X_2}\otimes\id_{E(X_3)}=
   d^E_{X_1, X_2\otimes X_3}\mcirc \id_{E(X_1)}\otimes d^E_{X_2,X_3}  \quad\forall X_1,X_2,X_3\in\2C
\]  
required from a tensor functor is immediate by associativity of $m$:
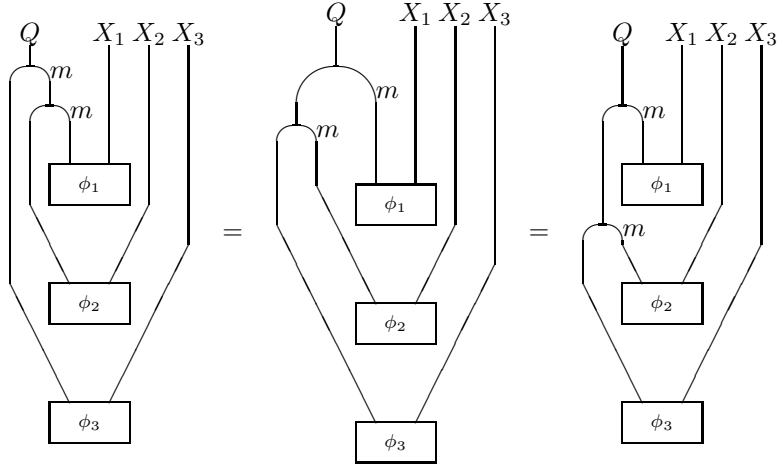
\[
\begin{tangle}
\hstep\object{Q}\Step\object{X_1}\step\object{X_2}\step\object{X_3}\\
\hh\hcd\obj{m}\step[1.5]\id\step\id\step\id\\
\hh\id\hstep\hcd\obj{m}\step\id\step\id\step\id\\
\hh\id\hstep\id\step\id\step\id\step\id\step\id\\
\hh\id\hstep\id\step\frabox{\phi_1}\step\id\step\id\\
\hh\id\hstep\d\Step\dd\step\id\\
\hh\id\step\d\step\dd\step\dd\\
\hh\d\step\frabox{\phi_2}\step\dd\\
\hstep\d\step\dd\\
\hh\step[1.5]\frabox{\phi_3}\\
\end{tangle}
\quad = \quad
\begin{tangle}
\step[1.5]\object{Q}\Step\object{X_1}\step\object{X_2}\step\object{X_3}\\
\hstep\cd\obj{m}\step\id\step\id\step\id\\
\hh\cd\obj{m}\step[1.5]\id\step\id\step\id\step\id\\
\hh\id\step\id\step[1.5]\id\step\id\step\id\step\id\\
\hh\id\step\d\step\frabox{\phi_1}\step\id\step\id\\
\hh\d\step\d\Step\dd\step\id\\
\hh\hstep\d\step\d\step\dd\step\dd\\
\hh\step\d\step\frabox{\phi_2}\step\dd\\
\step\hstep\d\step\dd\\
\hh\step\step[1.5]\frabox{\phi_3}\\
\end{tangle}
\quad = \quad
\begin{tangle}
\step\object{Q}\step[1.5]\object{X_1}\step\object{X_2}\step\object{X_3}\\
\hstep\hcd\obj{m}\step\id\step\id\step\id\\
\hh\hstep\id\step\id\step\id\step\id\step\id\\
\hh\hstep\id\step\frabox{\phi_1}\step\id\step\id\\
\hh\cd\obj{m}\Step\dd\step\id\\
\hh\id\step\d\step\dd\step\dd\\
\hh\d\step\frabox{\phi_2}\step\dd\\
\hstep\d\step\dd\\
\hh\step[1.5]\frabox{\phi_3}\\
\end{tangle}
\]
That $(E, (d_{X,Y}),e)$ satisfies the unit axioms is almost obvious. The first condition follows by 
\[ d_{X,\11}(\id_{E(X)}\otimes
   e)\phi=d_{X,\11}(\phi\otimes\eta)=m\otimes\id_X\mcirc\id_Q\otimes\phi\mcirc\eta=\phi, \] 
and the second is shown analogously.

So far we have shown that $E$ is a weak tensor functor for which 
$e: \11_{\mathrm{Vect}_\7F}\rarr E(\11 )$ is
an isomorphism. In order to conclude that $E$ is a (strong) tensor functor it remains to show that the
morphisms $d^E_{X,Y}$ are isomorphisms. Let $X,Y\in\2C$. We consider the bilinear map 
\bean \gamma_{X,Y}: && \Hom_{Q-\Mod}(Q,Q\otimes X)\boxtimes \Hom_{Q-\Mod}(Q,Q\otimes Y)\rarr
  \Hom_{Q-\Mod}(Q,Q\otimes X\otimes Y), \\
  && s\boxtimes t\mapsto s\otimes\id_Y\mcirc t, 
\eean
and we write $\boxtimes$ rather than $\otimes_\7F$ for the tensor product of $\Vect_\7F$ in order 
to avoid confusion with the tensor product in $Q-\Mod$. By 2. we have Q-module morphisms 
$s_i: Q\rarr Q\otimes X, s_i': Q\otimes X\rarr Q$ for  
$i=1,\ldots,n(X)$ satisfying $s_i'\circ s_j=\delta_{ij}\id_Q$,
and $\sum_i s_i\circ s_i'=\id_{Q\otimes X}$, and similar morphisms $t_i, t_i', \ i=1,\ldots,n(Y)$
with $X$ replaced by $Y$. 
Then the $\gamma_{ij}=\gamma_{X,Y}(s_i\otimes t_j)$ are linearly independent because 
$\gamma'_{i'j'}\circ\gamma_{ij}=\delta_{i'i}\delta_{j'j}\id_Q$, where 
$\gamma'_{i'j'}=t'_j \mcirc s'_i\otimes\id_Y$. Bijectivity of $\gamma_{X,Y}$ follows now from the
fact that both the domain and codomain of $\gamma_{X,Y}$ have dimension $n(X)n(Y)$.

For any $X\in\2C$ we have a $Q$-module $(Q\otimes X, m\otimes\id_X)$. If $(Q,m,\eta)$ is a monoid in
the tensor category then it is straightforward to check that the following maps are
inverses of each other: 
\bean \delta_X: && \Hom_{Q-\Mod}(Q,Q\otimes X)\rarr\Hom(\11,Q\otimes X), 
   \quad\quad s\mapsto s\circ\eta, \\
   \delta^{-1}_X: && \Hom(\11,Q\otimes X)\rarr\Hom_{Q-\Mod}(Q,Q\otimes X), \quad\quad 
   \tilde{s}\mapsto m\otimes\id_X \mcirc\id_Q\otimes\tilde{s}.
\eean
But
\[ d^E_{X,Y}= \delta_{X\otimes Y}\mcirc \gamma_{X,Y} \mcirc \delta_X^{-1} \boxtimes \delta_Y^{-1}, \]
which shows that $d^E_{X,Y}$ is an isomorphism for every $X,Y\in\2C$.
\qed

\brem
From the assumptions it follows that $Q\cong\oplus_i n(\ol{X_i})X_i$. Such an object $Q$ cannot
exist in $\2C$ if $\2C$ has infinitely many isomorphism classes of irreducible objects. This is the
reason why we consider monoids living in a larger category $\hat{\2C}$.
\erem

The previous considerations being valid over any field $\7F$, 
we now turn to $*$-categories where $\7F=\7C$.

\bprop \label{p-embed*}
Let $\2C$ be a semisimple tensor $*$-category and let $(Q,m,\eta)$ be a monoid in $\hat{\2C}$
satisfying the conditions of Proposition \ref{p-embed} and in addition: 
\begin{enumerate}
\setcounter{enumi}{2}
\item $Q\in\2C_*$.
\item For every $s\in\Hom_{Q-\Mod}(Q,Q\otimes X)$ we have $s^*\in\Hom_{Q-\Mod}(Q\otimes X,Q)$.
\end{enumerate}
Then the functor $E$ defined in Proposition \ref{p-embed} is $*$-preserving w.r.t.\ the scalar
products on $E(X)$ given by $(\phi,\psi)\id_\11=\psi^*\circ\phi$, and the isomorphisms $d_{X,Y}$ are
unitary for all $X,Y\in\2C$. 
\eprop
\prf Clearly the inner products are positive definite, thus the $E(X)$ Hilbert spaces. 
Let $s: X\rarr Y, \ \phi\in\Hom(\11,Q\otimes X)$ and $\psi\in\Hom(\11,Q\otimes Y)$. Then
\[ (E(s)\phi,\psi)=\psi^*\circ\id_Q\otimes s\circ\phi=(\phi^*\circ\id_Q\otimes s^*\circ\psi)^*=
   \ol{(E(s^*)\psi,\phi)}=(\phi,E(s^*)\psi). \]
Thus $E(s^*)=E(s)^*$, so $E$ is a $*$-preserving functor. 

By assumption 2 we have the isomorphism $Q\otimes X\cong n(X)Q$ in the category $Q-\Mod$, to wit 
there exist $s_i\in\Hom_{Q-\Mod}(Q,Q\otimes X),\ t_i\in\Hom_{Q-\Mod}(Q\otimes X,Q), i=1,\ldots,n(X)$
satisfying $t_i\circ s_j=\delta_{ij}\id_Q$ and $\sum_i s_i\circ t_i=\id_{Q\otimes X}$. Now
4. implies that we can choose the $s_i, t_i$ such that $t_i=s_i^*$. 
We must show that $d^E_{X,Y}:E(X)\otimes E(Y)\rightarrow E(X\otimes Y)$ is unitary
for every $X,Y\in\2C$. Since we already know that it is an isomorphism, it suffices to show that it
is an isometry. Since $\delta_X$ and $\delta_Y$ are isomorphisms, we need only show that
\[(d^E_{X,Y}(\delta_X(s_i )\otimes\delta_Y (s_j ))
,d^E_{X,Y}(\delta_X(s_{i'})\otimes\delta_Y (s_{j'})))_{E(X\otimes Y)}
=(\delta_X(s_i ),\delta_X (s_{i'} ))_{E(X)}(\delta_Y(s_j ),\delta_Y (s_{j'} ))_{E(Y)}\]
for all $i,i',j,j'$. But definition of the inner products, the R.H.S. equals
\[(\eta^*\circ s_{i'}^*\circ s_i\circ\eta)(\eta^*\circ s_{j'}^*\circ s_j\circ\eta )
=\delta_{ii'}\delta_{jj'}(\eta^*\circ\eta)^2 =\delta_{ii'}\delta_{jj'} ,\]
whereas the L.H.S. equals
\begin{eqnarray*} \lefteqn{ 
  (\delta_{X\otimes Y}\circ\gamma_{X,Y}(s_i\boxtimes s_j )
  ,\delta_{X\otimes Y}\circ\gamma_{X,Y}(s_{i'}\boxtimes s_{j'}))_{E(X\otimes Y)} } \\
 &&=(\delta_{X\otimes Y}\circ(s_i\otimes \id_Y )\circ s_j 
,\delta_{X\otimes Y}\circ(s_{i'}\otimes \id_Y )\circ s_{j'})_{E(X\otimes Y)} \\
 && =((s_i\otimes \id_Y )\circ s_j\circ\eta 
,(s_{i'}\otimes \id_Y )\circ s_{j'}\circ\eta )_{E(X\otimes Y)} \\
 && =\eta^*\circ s_{j'}^*\circ (s_{i'}^*\otimes\id_Y )\circ (s_i\otimes\id_Y )\circ s_j\circ\eta
  =\delta_{ii'}\delta_{jj'}\eta^*\circ\eta  =\delta_{ii'}\delta_{jj'}, 
\end{eqnarray*}
as desired.
\qed

\brem In the situation where $\2C=\Rep_f(A,\D)$ for a discrete AQG, we have seen that
$\hat{\2C}\simeq\Rep(A,\D)$ and $\2C_*\simeq\Rep_*(A,\D)$. The regular monoid
$(\pi_l,\tilde{m},\tilde{\eta})$ satisfies all assumptions of Proposition \ref{p-embed*}: As to
assumption 3, recall from Proposition \ref{p-disc} that $\pi_l\in\Rep_*(A,\D)$. 
Assumption 4 follows from unitarity of the isomorphism 
$V_\theta: \pi_l\times I_\theta\rarr\pi_l\times\theta$ and the fact that $V_\theta$ and $V_\theta^*$
are morphisms of $\pi_l$-modules.
\erem

\blemma Let $\2C$ be as in Proposition \ref{p-embed}. Let $Q\in\hat{\2C}$ be a direct sum of
irreducible objects in $\2C$ with finite multiplicities, where $\11$ appears with multiplicity
one. Consider the functor 
$\2C\rarr\Vect_\7F$ defined by $E(X)=\Hom_{\hat{\2C}}(\11,Q\otimes X)$. Then the map 
$a: \End\,Q\rarr\Nat\,E,\ \ s\mapsto (a_X(s))$ with $a_X(s)=s\otimes\id_X\in\End E(X)$, is an
isomorphism. It restricts to an isomorphism $\Aut\,Q\rarr\Aut\,E$.
\elemma

\prf That $(a_X(s))$ is a natural transformation from $E$ to itself is obvious. Injectivity follows
from $a_\11(s)=s\otimes\id_\11=s$. The fact $Q\cong\oplus_i n_i X_i$, where $i$ runs through $I$
and $n_i\in\7Z_+$, implies $\End\,Q\cong\prod_i M_{n_i}(\7F)$. On the other hand, by semisimplicity
of $\2C$ we have $\Nat\,E\cong \prod_i \End\,E(X_i)$, cf.\ e.g.\ \cite{MRT}. Now it is easy to see
that the composition of the latter two isomorphisms with the map $a: \End\,Q\rarr\Nat\,E$ preserves
the factors in the respective direct products. Then surjectivity follows from 
$\dim E(X_i)=n_i$.
\qed

\blemma \label{l-nat}
Let $\2C$ and the monoids $(Q,m,\eta)$ and $(Q',m',\eta')$ be as in Proposition \ref{p-embed}. 
Assume in addition that $\2C$ has 
duals and that $Q, Q'$ are direct sums of irreducibles in $\2C$ with finite multiplicities.  Let 
$E, E': \2C\rarr\Vect_\7F$ be the ensuing embedding functors. Then there is a bijection
between monoidal natural isomorphisms $b: E\rarr E'$ and isomorphisms $s: Q\rarr Q'$ of monoids.
\elemma

\prf One direction is easy: If $s: Q\rarr Q'$ is an isomorphism such that 
$s\circ m=m'\circ s\otimes s$ and $\eta=\eta'\circ s$, then we define 
$a_X(s): E(X)\rarr E'(X)$ by $a_X(s)\phi=s\otimes\id_X\mcirc\phi\in E'(X)$ for 
$\phi\in E(X)$. The family $(a_X)$ obviously is a natural isomorphism of $E$ and $E'$, and that it
is monoidal, i.e.\ satisfies $d^{E'}_{X,Y}\circ a_X\otimes a_Y=a_{X\otimes Y}\circ d^E_{X,Y}$ for
all $X,Y$, is obvious by the definition of $d^E, d^{E'}$ and the fact that $s$ is an isomorphism of 
monoids. 

As to the converse, the existence of a monoidal natural isomorphism $b: E\rarr E'$ implies
$\dim\Hom(\11,Q\otimes X)=\dim\Hom(\11,Q'\otimes X)$ for $X\in\2C$. By duality we have
$\dim\Hom(X_i,Q)=\dim\Hom(X_i,Q')$ for all irreducible $X_i\in\2C$, which implies that 
$Q$ and $Q'$ are isomorphic. Fix an arbitrary isomorphism $s: Q'\rarr Q$ and consider the
monoid $(Q,m'',\eta'')$ where $m''= s\circ m'\circ s^{-1}\otimes s^{-1}$ and 
$\eta''=s\circ\eta'$. Let $E''$ be the embedding functor corresponding to $(Q,m'',\eta'')$.
By construction, $(Q,m'',\eta'')\cong (Q',m',\eta')$, and by the preceding
considerations we have the monoidal natural isomorphism $a(s)=(a_X(s)): E'\rarr E''$.
If $b: E\rarr E'$ is a monoidal natural isomorphism, then the composition 
$c=a(s)\circ b: E\rarr E''$ is monoidal, and there exists  
$t\in\Aut\,Q$ such that $c=c(t)$. Since $E$ and $E''$ coincide as functors, the condition
$c_{X\otimes Y}\circ d^E_{X,Y}=d^{E''}_{X,Y}\circ c_X\otimes c_Y$ is equivalent to
\[ (t\circ m)\otimes\id_{X\otimes Y}\mcirc\id_Q\otimes\phi\otimes\id_Y
   \mcirc\psi = (m''\circ t\otimes t)\otimes\id_{X\otimes Y}\mcirc
  \id_Q\otimes\phi\otimes\id_Y\mcirc\psi \]
for $X,Y\in\2C$ and $\phi\in E(X),\, \psi\in E(Y)$. Since $\2C$ has duals,
this means that
\[ m''\mcirc t\otimes t\mcirc u\otimes v=t\mcirc m\mcirc u\otimes v \]
for $X,Y\in\2C$ and $u: X\rarr Q,\ v: Y\rarr Q$. But $Q$ is a direct sum of simple objects in
$\2C$, so we can cancel $u\otimes v$ and conclude $m''\mcirc t\otimes t=t\mcirc m$. The equality
$\eta''\mcirc t=\eta$ is proven in a similar fashion using the morphisms 
$e^E: \7F\rarr E(\11)$ and $e^{E''}: \7F\rarr E''(\11)$. Thus we have an isomorphism 
$t: (Q,m,\eta)\rarr(Q,m'',\eta'')$ of monoids and composing with the isomorphism 
$s^{-1}: (Q,m'',\eta'')\rarr (Q',m',\eta')$ implies the claim. Clearly this gives us a
bijection between isomorphisms of embedding functors and of monoids, respectively.
\qed

\brem Having assumed throughout that the tensor category $\2C$ is strict, we now comment briefly 
on the non-strict case. If a tensor category $\2C$ has a non-trivial associativity constraint 
$$\alpha_{X,Y,Z}:(X\otimes Y)\otimes Z\rarr X\otimes(Y\otimes Z),$$ 
the definition of a monoid in $\2C$ is changed in an obvious way: 
The associativity condition becomes
$$m\circ(m\otimes\id_Q)=m\circ(\id_Q\otimes m)\circ\alpha_{Q,Q,Q},$$
and the first equation in Definition \ref{d-module} relating elements in 
$\Hom((Q\otimes Q)\otimes X,X)$ becomes
$$\mu\mcirc m\otimes\id_X=\mu\mcirc \id_Q\otimes\mu\mcirc\alpha_{Q,Q,X}.$$ 
It then remains true that an absorbing monoid gives rise to an embedding functor, 
but we omit the proofs. 
\erem


\subsection{Main Result}
Given a discrete AQG it is occasionally convenient to consider an abstract tensor
$*$-category $\Rep^{abs}_f(A,\D)$ together with an embedding functor $E$, rather than the concrete
category $\Rep_f(A,\D)$ and the forgetful functor $K$. 

\blemma \label{l-circ}
Let $(A,\D)$ be a discrete AQG and write $\2C=\Rep^{abs}_f(A,\D)$. Let $E: \2C\rarr\2H$ be the
obvious embedding functor. Let $(\pi_l,\tilde{m},\tilde{\eta})$ be the regular monoid in
$\hat{\2C}\simeq\Rep(A,\D)$ and $E': \2C\rarr\2H$ the embedding functor that it gives rise to by
Proposition \ref{p-embed*}. Then there exists a unitary equivalence $u: E\rarr E'$ of tensor functors.
\elemma

\prf
For $X\in\2C$ we have $(E(X),\pi_X )\in\Rep_f (A,\Delta )$, and let us write $V_X$ instead of 
$V_{\pi_X}$. For $\phi\in E(X)$ define $u_X\phi\in A\otimes E(X)$ by
$u_X\phi =V_X (I_0\otimes\phi )$. Then
\begin{multline*}  (\pi_l\times\pi_X )(a)u_X\phi   =(\pi_l\times\pi_X )(a)V_X (I_0\otimes\phi )
   =V_X (\pi_l\times I_{\pi_X})(a)(I_0\otimes\phi )
   =V_X (\pi_l\otimes\varepsilon )\Delta (a)(I_0\otimes\phi) \\
  =V_X (\pi_l (a)I_0\otimes\phi) =V_X (\varepsilon (a)I_0\otimes\phi )
   =\varepsilon (a)V_X (I_0\otimes\phi )=\varepsilon (a)u_X\phi,
\end{multline*}
thus $u_X\phi\in\Hom (\varepsilon ,\pi_l\times\pi_X )$.  
In order to show that $(u_X )$ is a natural
transformation, we consider $s:X\rightarrow X'$ and compute
\[u_{X'}E(s)\phi =V_{X'}(I_0\otimes s\phi )=V_{X'}(1\otimes s)(I_0\otimes\phi )
=(1\otimes s)V_X (I_0\otimes\phi )=(1\otimes s)u_X\phi =E' (s)u_X\phi ,\]
where we have used Proposition \ref{p-functorial}. 
Since $V_X$ is invertible, the map $\phi\mapsto u_X\phi$ is injective and therefore bijective by
equality of the dimensions.

Thus $(u_X )$ is a natural isomorphism.
It remains to show that it is monoidal, i.e.
\[d^{E'}_{X,X'}\mcirc u_X\otimes u_{X'}=u_{X\otimes X'} \]
for $X,X'\in\2C$. Here we have as usual identified the vector spaces $E(X)\otimes E(X' )$
and $E(X\otimes X' )$. Let $\phi\in E(X)$ and $\phi'\in E(X' )$.
Then 
\[u_{X\otimes X'}(\phi\otimes\phi')=V_{X\otimes X'}(I_0\otimes\phi\otimes\phi' ),\]
whereas 
\begin{eqnarray*} d^{E'}_{X,X'}\circ (u_X\otimes u_{X'})(\phi\otimes\phi ')
   &=& (\tilde{m}\otimes\id_X\otimes\id_{X'})\circ (\id_A\otimes u_X\phi\otimes\id_{X'})\circ
    u_{X'}\phi' \\
   &=&(\tilde{m}\otimes\iota\otimes\iota )
((V_{X'})_{14}(V_X)_{23}(1\otimes 1\otimes \phi\otimes\phi' )).
\end{eqnarray*}
Thus we must show that
\[V_{X\otimes X'}(I_0\otimes\phi\otimes\phi' ) 
=(\tilde{m}\otimes\iota\otimes\iota )
((V_{X'})_{14}(V_X)_{23}(1\otimes 1\otimes \phi\otimes\phi' )).\]
By non-degeneracy of $\pi_X$ and $\pi_{X'}$ we may assume $\phi =\pi_X (a)v$ 
and $\phi' =\pi_{X'} (b)v'$, for $a,b\in A$ and $v\in E(X)$ and $v'\in E(X' )$.
By the definition of $V_X$, $V_{X'}$ and $V_{X\otimes X'}$ it thus suffices to show that
\[(\tilde{m}\otimes\iota\otimes\iota )
(\Delta (I_0 )_{14}\Delta (I_0 )_{23}(1\otimes 1\otimes a\otimes b ))
=(\Delta\otimes\iota )\Delta (I_0 )(1\otimes a\otimes b) \]
for $a,b\in A$. Write  
$\Delta (I_0 )(1\otimes b)=\sum_i a_i\otimes b_i$ for $a_i ,b_i\in A$.
Then by Proposition \ref{p-qf} and $\tilde{m}(\id\otimes\tilde{\eta})=\id_{\pi_l}$, we get
\begin{eqnarray*} \lefteqn{  (\tilde{m}\otimes\iota\otimes\iota)
  (\Delta (I_0 )_{14}\Delta (I_0 )_{23}(1\otimes 1\otimes a\otimes b )) 
  = \sum_i (\tilde{m}\otimes\iota\otimes\iota )(a_i\otimes\Delta (I_0 )(1\otimes a)\otimes b_i) } 
    \quad\quad\\
  &&  =\sum_i (\tilde{m}\otimes\iota )(a_i\otimes\Delta (I_0 )(1\otimes a))\otimes b_i 
   =\sum_i \Delta\tilde{m}(a_i\otimes I_0 )(1\otimes a)\otimes b_i \\
  && =\sum_i \Delta(a_i )(1\otimes a)\otimes b_i 
   =\sum_i (\Delta\otimes\iota )(a_i\otimes b_i )(1\otimes a\otimes 1) \\
  && =(\Delta\otimes\iota )(\Delta (I_0 )(1\otimes b))(1\otimes a\otimes 1)
  =(\Delta\otimes\iota )\Delta (I_0 )(1\otimes a\otimes b),
\end{eqnarray*}
as desired. 
\qed

At this stage we need to recall the generalized Tannaka theorem for discrete AQG, as proven in
\cite{MRT}. 

\btheor
Let $\2C$ be a semisimple tensor $*$-category and let $E$ be a an embedding functor. Then there
exists a discrete AQG $(A,\D)$ and an equivalence $F: \2C\rarr\Rep_f(A,\D)$ of tensor
$*$-categories, such that $K\circ F=E$, where $K: \Rep_f(A,\D)\rarr\2H$ is the forgetful functor. 
\etheor

We are now in a position to state our main result which describes
the precise relationship between embedding functors, absorbing monoids and discrete AQG.

\btheor \label{c-equiv}
\ 1. Let $\2C$ be a tensor $*$-category with conjugates and $\End\,\11\cong\7C$ and let $E:
\2C\rarr\2H$ be an embedding functor. Let $(A,\D)$ be the discrete AQG and $F: \2C\rarr\Rep_f(A,\D)$
the monoidal equivalence provided by the generalized Tannaka theorem. 
Let $(\pi_l,\tilde{m},\tilde{\eta})$
be the regular monoid in $\Rep(A,\D)$ and $E': \2C\rarr\2H$ the embedding functor that it gives rise
to. Then  $E$ and $E'$ are naturally unitarily equivalent as tensor functors.

2. Let $(A,\D)$ be a discrete AQG and $(\pi_l,\tilde{m},\tilde{\eta})$ the regular monoid in 
$\Rep(A,\D)$. Let $E: \Rep^{abs}_f(A,\D)\rarr\2H$ be the embedding functor obtained from the latter
via Proposition \ref{p-embed} and $(A',\D')$ the discrete AQG given by the generalized Tannaka
theorem. Then $(A,\D)$ and $(A',\D')$ are isomorphic.  

3. Let $\2C$ be a tensor $*$-category with conjugates and $\End\,\11\cong\7C$ and let $(Q,m,\eta)$
be a monoid in $\hat{\2C}$ satisfying the assumptions in Proposition \ref{p-embed*}. Let $E$ be the
resulting embedding functor and $(A,\D)$ and $F$ as in 2. Then the image $(Q',m',\eta')$
of the regular monoid $(\pi_l, \tilde{m},\tilde{\eta})$ under the equivalence
$\Rep(A,\D)\rarr\hat{\2C}$ is isomorphic to $(Q,m,\eta)$.
\etheor

\prf 1. Consider the equivalence
$F: \2C\rarr\Rep_f(A,\D)$ satisfying $K\circ F=E$ provided by the generalized Tannaka theorem. 
Then the claim is just a reformulation of Lemma \ref{l-circ}.

2. Let $\2C=\Rep_f^{abs}(A,\D)$ with the canonical embedding
functor $E: \2C\rarr\2H$. Obviously, $(A,\D)$ is isomorphic to the AQG given by the generalized
Tannaka theorem from the pair $(\2C,E)$. Now the claim follows from Lemma \ref{l-circ} and the
fact \cite[Proposition 5.28]{MRT} that isomorphic embedding functors give rise to isomorphic
discrete AQG. 

3. Given $\2C$ and the monoid $(Q,m,\eta)$ in $\hat{\2C}$, we obtain an embedding functor 
$E: \2C\rarr\2H$ by Proposition \ref{p-embed*}. On the other hand, going from $(\2C,E)$ to an AQG, 
then to the regular monoid in $\Rep(A,\D)\simeq\hat{\2C}$ and, finally, from the latter to
the embedding functor $E': \2C\rarr\2H$, Lemma \ref{l-circ} again implies
$E\stackrel{\otimes}{\cong}E'$. Thus the monoids $(Q,m,\eta)$ and $(Q',m',\eta')$ in $\hat{\2C}$
give rise to equivalent embedding functors and are therefore isomorphic by Lemma \ref{l-nat}.
\qed

\brem 1. The preceding result can be formalized more conceptually as follows. Let $Disc$ be the
category of discrete AQG with isomorphisms as arrows. Let $Emb$ be the category of  pairs $(\2C,E)$
where $\2C$ is a semisimple $\7F$-linear tensor category with duals and $\End\,\11\cong\7F$ and 
$E: \2C\rarr\Vect_\7F$ is a faithful $\7F$-linear tensor functor. The arrows in $\mathrm{Emb}$ are
equivalences $F: \2C\rarr\2C'$ such that $E'\circ F=E$. Finally, let $\mathrm{Mon}$ be the category
of pairs $(\2C,(Q,m,\eta))$, where $\2C$ is a semisimple $\7F$-linear tensor category with duals and 
$\End\,\11\cong\7F$ and $(Q,m,\eta)$ is a monoid in $\hat{\2C}$ satisfying the assumptions of
Proposition \ref{p-embed}. Here the arrows are equivalences $F: \2C\rarr\2C'$ such that
$F((Q,m,\eta))$ is isomorphic to $(Q',m',\eta')$ in $\2C'$. Then the various constructions
considered so far give rise to the equivalences 
$\mathrm{Mon}\simeq\mathrm{Emb}\simeq\mathrm{Disc}^{op}$, where $\mathrm{Disc}^{op}$
is the opposite category of $\mathrm{Disc}$. More precisely, every circle in the triangle with
corners $\mathrm{Mon},\mathrm{Emb}, \mathrm{Disc}^{op}$ obtained as composition of these functors is
naturally isomorphic to the identity functor. 

2. The preceding theorem remains valid if one replaces tensor $*$-categories with conjugates
by semisimple $\7F$-linear tensor categories with duals, and discrete AQG by regular multiplier Hopf
algebras with left invariant functionals. The arguments are essentially unchanged, provided one
appeals to the version of the generalized Tannaka theorem stated in \cite[Section 5.4]{MRT}.
\erem


\subsection{Dimension Functions vs.\ Absorbing Objects}
\bdefin A dimension function on a $C^*$-tensor category $\2C$ with conjugates is a map 
$n:\obj\,\2C\rarr\7R_+$ such that $n(X\oplus Y)=n(X)+n(Y)$ and $n(X\otimes Y)=n(X)n(Y)$ and 
$n(X)=n(\ol{X})$. 
\edefin

\brem
Note that a dimension function automatically satisfies $n(\11)=1$. Every $C^*$-tensor category 
$\2C$ with conjugates comes with a distinguished dimension function, the intrinsic dimension, cf.\
\cite{LR}. The representation categories associated with $q$-deformations of simple Lie groups show
that the intrinsic dimension need not be integer valued, cf.\ \cite{RT}.
On the other hand, an embedding functor
$E:\2C\rarr\2H$ gives rise to an integer valued dimension function by $n(X)=\dim E(X)$. This also
shows that one and the same category can have a dimension function which is integer valued and one
which is not. We remark further that $C^*$-tensor categories having only finitely many irreducible
objects admit only one dimension function, namely the intrinsic one, as can be shown using
Perron-Frobenius theory. Furthermore, every embedding functor must preserve dimensions whenever
$\2C$ is amenable, which in particular holds when $\2C$ admits a unitary braiding, cf.\
\cite{LR}. Thus if the intrinsic dimension of $\2C$ is not integer valued and $\2C$ is finite or has
a unitary braiding, an embedding functor cannot exist.
\erem

Assuming the existence of an integer valued dimension function 
we arrive at the following partial converse of Proposition \ref{p-embed}. 

\bprop \label{l-absor}
Let $\2C$ be a semisimple $\7F$-linear tensor category with two-sided duals and integer valued
dimension function $n$. Let $n_i=n(X_i)$ for $i\in I_\2C$ and consider the direct sum  
\[ Q=\bigoplus_{i\in I_\2C} n_{\ol{\imath}}\,X_i \]
in $\hat{\2C}$. Then $Q\otimes X\cong X\otimes Q\cong n(X) Q$
for all $X\in\2C$. 

Conversely, assume $Q\in\hat{\2C}$ is a direct sum of irreducible objects of $\2C$ and that 
$Q\otimes X\cong n(X) Q$ with $n(X)\in\7N$ for $X\in\2C$. Then 
\[ Q \ \cong\ N\,\bigoplus_{i\in I_\2C} n_{\ol{\imath}}\, X_i, \]
where $N=\dim\Hom(\11,Q)$. 
If $N<\infty$ then $n: \obj\,\2C\rarr\7N$ is additive and multiplicative.
If $I_\2C$ is a finite set then $n(X)=d(X)$ for all $X\in\2C$; thus in this case an absorbing object
exists iff all intrinsic dimensions are integers. 
\eprop

\prf By $\bigoplus_{i\in I_\2C} n_{\ol{\imath}} X_i$ we mean the filtered inductive limit over
partial finite direct sums, which defines an object of $\hat{\2C}$ unique up to isomorphism.
Let $j\in I_\2C$. We compute
\[ Q\otimes X_j\cong\bigoplus_{i\in I_\2C} n_{\ol{\imath}}\,X_i\otimes X_j 
  \cong\bigoplus_{i\in I_\2C}n_{\ol{\imath}}\bigoplus_{k\in I_\2C}
   N_{ij}^k X_k \cong \bigoplus_{k\in I_\2C} (\sum_{i\in I_\2C} N_{ij}^k n_{\ol{\imath}} ) \, X_k. 
\]
Using standard properties of the coefficients, cf.\ e.g.\ \cite{MRT}, we calculate
\[ \sum_{i\in I_\2C} N_{ij}^k n_{\ol{\imath}}= \sum_{i\in I_\2C}
  N_{j\ol{k}}^{\ol{\imath}}n_{\ol{\imath}}=  n_j n_{\ol{k}}, \]
and therefore $Q\otimes X_j\cong n_j\bigoplus_{k\in I_\2C} n_{\ol{k}} X_k\cong n_jQ$. For a reducible
object $X$ the claim now follows by semisimplicity. The argument for $X\otimes Q$ is similar.

As to the converse, for irreducible $X\in\2C$ we compute
\[ \dim\Hom(X,Q)=\dim\Hom(\11,Q\otimes\ol{X})=\dim\Hom(\11,Q\otimes \11^{\oplus n(\ol{X})})
    =n(\ol{X})\dim\Hom(\11,Q)=n(\ol{X})N. \]
Since $Q$ is a direct sum of irreducibles in $\2C$, we thus have
\[ Q\cong N\bigoplus_{i\in I_\2C} n_{\ol{\imath}} X_i, \]
and the claim follows. Assume now that $N<\infty$. Then we find
\begin{eqnarray*} n(X\otimes Y)N &=& \dim\Hom(\ol{X\otimes Y},Q) \ =\
   \dim\Hom(\ol{Y}\otimes\ol{X},Q) \\ 
  &=& \dim\Hom(\ol{Y},Q\otimes X)\ =\ n(X)\dim\Hom(\ol{Y},Q)\ =\ n(Y)n(X)N 
\end{eqnarray*}
and thus $n(X)n(Y)=n(X\otimes Y)$ for $X,Y\in\2C$.

If $\2C$ is finite, it is well known that the intrinsic dimension function is the only additive and
multiplicative function on $\obj\,\2C$.
\qed

\brem
1. Note that an additive and multiplicative function on $\obj\,\2C$ 
determines and is determined by a function $n': I_\2C\rarr\7N$ which satisfies 
$\sum_{k\in I_\2C} N_{ij}^k n'_k=n'_i n'_j$ for all $i,j\in I_\2C$.

2. It is important to note that the existence of an integer valued dimension function does not 
obviously imply the existence of a monoid structure on the absorbing object $Q$. By our earlier
constructions, an embedding functor gives rise to a quantum group, and therefore to the regular 
monoid in $\hat{\2C}$. (One can also construct the latter directly from the embedding functor, 
but we refrain from giving the details.) Since any dimension function $n$ satisfies $n(\11)=1$, 
we have $\dim\Hom(\11,Q)=1$, thus there exists a morphism $\eta:\11\rarr Q$ that is unique up to a
scalar. But the main issue clearly is constructing an associative morphism $m:Q\otimes Q\rarr Q$
such that $(Q,m,\eta)$ is a monoid. This is a difficult cohomological problem.

3. Another approach for constructing an absorbing monoid might be to generalize Deligne's
proof to the braided case. However, as our earlier mentioned counter examples show,
assuming just the existence of a braiding does not suffice. 
\erem



\begin{thebibliography}{30}

\bibitem{SGA4} M. Artin, A. Grothendieck \& J. L. Verdier: Th\'{e}orie des topos et cohomologie
\'{e}tale des sch\'{e}mas. Tome 1: Th\'{e}orie des topos. (SGA 4). Springer, 1972.

\bibitem{bichon} J. Bichon: Trivialisations dans les cat\'egories tannakiennes. 
  Cah. Topol. Geom. Cat\'eg. {\bf 39}, 243--270 (1998).

\bibitem{del} P. Deligne: Cat\'{e}gories tannakiennes. {\it In:} P. Cartier et al. (eds.):
{\it Grothendieck Festschrift}, vol. II, pp. 111--195. Birkh\"auser, 1991.

\bibitem{del2} P. Deligne: Le groupe fondamental de la droite projective moins trois points. 
{\it In:} Y. Ihara, K. Ribet \& J.-P. Serre (eds.): {\it Galois groups over} $\7Q$,
Math. Sci. Res. Inst. Publ. {\bf 16}, 79--297. Springer, 1989.

\bibitem{DM} P. Deligne \& J. S. Milne: Tannakian categories. 
Lecture notes in mathematics {\bf 900}, 101--228. Springer, 1982.

\bibitem{DR} S. Doplicher \& J. E. Roberts: A new duality theory for compact groups.
Invent. Math. {\bf 98}, 157--218 (1989).

\bibitem{DR2} S. Doplicher \& J. E. Roberts: Why there is a field algebra with compact gauge group
describing the superselection structure in particle physics.  
Commun. Math. Phys. {\bf 131}, 51--107 (1990).

\bibitem{dpr} S. Doplicher, C. Pinzari \& J. E. Roberts: An algebraic duality theory for 
multiplicative unitaries. Int. J. Math. {\bf 12}, 415--459 (2001). C. Pinzari \& J. E. Roberts: 
Regular objects, multiplicative unitaries and conjugation.  Internat. J. Math. {\bf 13}, 625-665 (2002).

\bibitem{ddz} B. Drabant, A. Van Daele \& Y. H. Zhang: Actions of multiplier Hopf algebras. 
Commun. Alg. {\bf 27}, 4117--4172 (1999).

\bibitem{FI} F. Fidaleo \& T. Isola: The canonical endomorphism for infinite index inclusions.
Z. Anal. Anw. {\bf 18}, 47-66 (1999).

\bibitem{gab} P. Gabriel: Des cat\'{e}gories ab\'{e}liennes. 
Bull. Soc. Math. France {\bf 90}, 323--448 (1962).

\bibitem{JS} A. Joyal \& R. Street: An introduction to Tannaka duality and quantum groups.
{\it In:} Proceedings of the conference on category theory, Como 1990. 
Lecture notes in mathematics {\bf 1488}, 413--492. Springer, 1991.

\bibitem{kas} C. Kassel: {\it Quantum Groups}. Springer, 1995.

\bibitem{Wenzl2} D. Kazhdan \& H. Wenzl: Reconstructing tensor categories. 
Adv. Sov. Math. {\bf 16}, 111--136 (1993). 

\bibitem{KD} J. Kustermans \& A. Van Daele: C$^{*}$-algebraic quantum groups arising from algebraic 
quantum groups. Int. J. Math. {\bf 8}, 1067--1139 (1997).

\bibitem{longo} R. Longo: A duality for Hopf algebras and for subfactors. I. 
Commun. Math. Phys. {\bf 159}, 133-150 (1994).

\bibitem{LR} R. Longo \& J. E. Roberts: A theory of dimension. 
K-Theory {\bf 11}, 103--159 (1997).

\bibitem{cwm} S. Mac Lane: {\it Categories for the Working Mathematician}. 2nd edition. Springer, 1997.

\bibitem{mue09} M. M\"uger: From subfactors to categories and topology I. Frobenius algebras in and
Morita equivalence of tensor categories. J. Pure Appl. Alg. {\bf 180}, 81--157 (2003).

\bibitem{MM} M. M\"uger: Abstract duality for symmetric tensor *-categories. Appendix to: Hans
Halvorson: Algebraic Quantum Field Theory. In: J. Butterfield and J. Earman (eds.): {\it
Handbook of the Philosophy of Physics}. Elsevier, 2006. ({\tt math-ph/0602036}).

\bibitem{MRT} M. M\"uger, J. E. Roberts \& L. Tuset: Representations of algebraic quantum groups 
and reconstruction theorems for tensor categories. Alg. Repres. Theor. {\bf 7}, 517--573 (2004).

\bibitem{RT} J. E. Roberts \& L. Tuset: On the equality of $q$-dimensions and intrinsic dimension. 
J. Pure Appl. Alg. {\bf 156}, 329--343 (2001).

\bibitem{Wenzl} H. Wenzl: $C^*$-Tensor categories from quantum groups. 
  J. Amer. Math. Soc. {\bf 11}, 261--282 (1998).

\bibitem{woro} S. L. Woronowicz: Tannaka-Krein duality for compact matrix pseudogroups. Twisted 
SU(N) groups. Invent. Math. {\bf 93}, 35--76 (1988). 

\bibitem{VD2} A. Van Daele: Multiplier Hopf algebras. 
Trans. Amer. Math. Soc. {\bf 342}, 917--932 (1994).

\bibitem{VD} A. Van Daele: An algebraic framework for group duality.
 Adv. Math. {\bf 140}, 323--366 (1998).

\bibitem{yama1} S. Yamagami: Group symmetry in tensor categories and duality for orbifolds. 
J. Pure Appl. Alg. {\bf 167}, 83-128 (2002).

\bibitem{yama2} S. Yamagami: $C^*$-tensor categories and free product bimodules. 
J. Funct. Anal. {\bf 197}, 323-346 (2003). 
\end{thebibliography}
\end{document}